\documentclass[11 pt, psamsfonts]{amsart}
\usepackage{amssymb}
\usepackage{amsmath}	
\usepackage{ bm}
\usepackage{ytableau}
\usepackage[normalem]{ulem}
\usepackage{verbatim}
\usepackage{array,calc,url,tabularx, epic,tikz}
\usepackage{comment}
\usetikzlibrary{arrows,positioning,decorations.pathmorphing,
  decorations.markings}
\usetikzlibrary{arrows.meta}

\def\ignore#1{\relax}

\def\g{\mathfrak g}
\def\h{\mathfrak h}

\def\wt{{\rm{wt}}}

\def\R{{\mathbb R}}
\def\Z{{\mathbb Z}}

\def\C{{\mathbb C}}

\def\la{\lambda}
\def\La{\Lambda}

\def\G{\mathcal G}
\def\Neg{\mathcal Neg}

\def\N{\mathbb N}
\def\A{\mathcal A}
\def\Ca{\mathcal C}

\def\P{{\bf P}}

\def\A{\mathcal A}

\def\a{{\bf a}}

\def\U{\mathcal{U}}

\def\ni{\noindent}

\def\ignore#1{\relax}

\def\om{\omega}

\def\eps{\epsilon}
\def\1{{\bf 1}}

\def\ep{\epsilon}

\def\End{{\rm End}}
\def\Hom{{\rm Hom}}
\def\ve{\varepsilon}
\def\eps{\varepsilon}

\calclayout
\renewenvironment{proof}[1][\proofname]{\par
  \normalfont
  \topsep6\p@\@plus6\p@ \trivlist
  \item[\hskip\labelsep\bfseries
    #1\@addpunct{.}]\ignorespaces
}{%
  \qed\endtrivlist
}

\def\th@plain{%
  \let\thmhead\thmhead@plain \let\swappedhead\swappedhead@plain
  \thm@preskip.5\baselineskip\@plus.2\baselineskip
                                    \@minus.2\baselineskip
  \thm@postskip\thm@preskip
  \itshape
\renewcommand{\labelenumi}{{(\alph{enumi})\quad}}
                        \renewcommand{\labelenumii}{{(\roman{enumii})\ }}
}
\def\th@definition{%
  \let\thmhead\thmhead@plain \let\swappedhead\swappedhead@plain
  \thm@preskip.5\baselineskip\@plus.2\baselineskip
                                    \@minus.2\baselineskip
  \thm@postskip\thm@preskip
  \upshape
}
\def\th@remark{%
  \thm@headfont{\itshape}% heading font bold
  \let\thmhead\thmhead@plain \let\swappedhead\swappedhead@plain
  \thm@preskip.5\baselineskip\@plus.2\baselineskip
                                    \@minus.2\baselineskip
  \thm@postskip\thm@preskip
  \upshape
}

{\theoremstyle{plain}
\newtheorem{theorem}{Theorem}[section]
}

{\theoremstyle{plain}
\newtheorem{proposition}[theorem]{Proposition}
}

{\theoremstyle{plain}
\newtheorem{corollary}[theorem]{Corollary}
}

{\theoremstyle{plain}
\newtheorem{lemma}[theorem]{Lemma}
}

{\theoremstyle{plain}

}

{\theoremstyle{definition}
\newtheorem{definition}[theorem]{Definition}
}

{\theoremstyle{definition}

}

{\theoremstyle{remark}
\newtheorem{remark}[theorem]{Remark}
}

{\theoremstyle{remark}

}

\numberwithin{equation}{section}

\renewcommand{\labelenumi}{{ \theenumi.}}
\renewcommand{\labelenumii}{{(\alph{enumii})}}

\def\ni{\noindent}

\def\la{\lambda}
\def\al{\alpha}

\def\ga{\gamma}

\def\choose #1 #2{\begin{pmatrix}#1\\#2\end{pmatrix}}

\begin{document}
%\frontmatter

\title[$G_2$ categories]
{Reconstruction of tensor categories of type $G_2$}

\author{Lilit Martirosyan }
\address{L.M. Department of Mathematics and Statistics\\ University of North Carolina\\ Wilmington \\North Carolina}
\email{martirosyanl@uncw.edu}

\author{Hans Wenzl}
\address{H.W. Department of Mathematics\\ University of California\\ San Diego,
California}

\email{hwenzl@ucsd.edu}

\begin{abstract} We prove that any ribbon tensor category $\Ca$ with the fusion rules of the compact group of type $G_2$ is equivalent to the representation category of the corresponding Drinfeld-Jimbo quantum group for $q^2=1$ or $q^2$ not a root of unity. We also prove an analogous result for the corresponding ribbon fusion tensor categories at roots of unity.
\end{abstract}
\maketitle
Given an index set for the equivalence classes of simple objects of a tensor category and a collection of possible tensor product rules for them, it is a natural question to ask how many tensor categories satisfy these fusion rules. If $G$ is a compact group, we say that a semisimple tensor category $\Ca$ is of type $G$ if its fusion semiring is isomorphic to that of Rep $G$. 
It was shown in \cite{KW}
that any category of type $SU(N)$ has to be equivalent to Rep $U_qsl_N$ for $q$ not a root of unity, whose associativity morphism could be twisted by a 3-cocycle of the center $\Z/N$ of $SU(N)$. Similar results were also obtained for the associated fusion tensor categories $SU(N)_k$ where now $q$ has to be a primitive $(k+N)$-th root of unity.
An analogous result was proved in \cite{TW} for groups $G=O(N)$ and $Sp(N)$, and in \cite{Cp} for $SO(N)$, assuming that the tensor categories are ribbon tensor categories. 

In this paper, we deal with ribbon tensor categories $\Ca$ of type $G(G_2)$ (or just of type $G_2$ for short), the compact group of Lie type $G_2$, and the associated finite fusion tensor categories of type $G_{2,k}$, see Section \ref{fusion:sec} for precise definitions. Our main result is as follows.
\vskip .2cm
\ni {\bf Main theorem} (see Theorems \ref{classification-theorem} and \ref{explicit classification}) (a) If $\Ca$ is a ribbon tensor category of type $G_2$, it is equivalent to $\U_q=$ Rep $U_q\g(G_2)$, with $q^2=1$ or $q^2$ not a root of unity.

(b) If $\Ca$ is a ribbon tensor category of type $G_{2,k}$ with $k\geq 1$, $k\neq 3,6$, it is equivalent to the quotient category $\bar{\U}_q$ of $\U_q$ with $q^2$ a primitive $(k+12)$-th root of unity.

(c) For $k=3$, the categories of type $G_{2,3}$ are
the two Fibonacci categories. For $k=6$, the categories of
type $G_{2,6}$ are the categories with the fusion rules of
$SO(3)_7$, classified as in Lemma \ref{lowtensor}.

(d) Let $\Ca(q^2)$ denote the categories in
parts (a) and (b). Then $\Ca(q^2)$ is equivalent to $\Ca(\tilde q^2)$ as a monoidal tensor category if and only if $\tilde q^2=q^{\pm 2}$, except for type $G_{2,9}$, where we have the additional equivalences $\Ca(q^{\pm 2})\sim\Ca(q^{\pm 8})\sim \Ca(q^{\pm 10})$ for $q^2$ a primitive 21st root of unity.

\vskip .2cm

A major ingredient in our proof is the complete classification of all simple representations of the quotient $K_4$ of the braid group algebra $\C B_4$ for which the standard generators are diagonalizable and satisfy a cubic equation. This appears in our companion paper \cite{MW3}. It allows us to
show that the $B_4$ representation on $\Hom(V_{\La_1+\La_2},V^{\otimes 4})$ must be the 8-dimensional canonical indecomposable representation of $K_4$; here $V$ is the object corresponding to the simple 7-dimensional representation of $G_2$. We deduce from this the eigenvalues of the braiding element $c_{V,V}$.

A second major tool is what we call path rigidity. Using a generalization of quantum Jucys-Murphy elements, see Theorem \ref{AB2theorem}, and the techniques in \cite{MW2}, we show that the canonical representations $\C B_n\to \End(V^{\otimes n})$ are surjective, and for type $G_2$ also uniquely determined by the eigenvalue $q^2$ of $c_{V,V}$ for $V_{2\La_1}\subset V^{\otimes 2}$. This can be extended to tensor categories of type $G_{2,k}$ with some effort, see Remark \ref{general case:rem}.

There is another approach to categories of type $G_2$, based on Kuperberg's spiders, see \cite{Ku2} and subsequent work in \cite{MPS} and \cite{MST}. Connecting our description via braids to the one with spiders allows us to give a second, complete proof of our main result avoiding some of the technicalities in \cite{MW2}. In particular, we can re-express our main result as saying that all categories of type $G_2$ or $G_{2,k}$ are equivalent to spider (or tri-valent) categories modulo negligible morphisms, see Theorem \ref{quotient:corollary}.

Here is our paper in more detail. We review basics of tensor categories and quantum groups of type $G_2$ in the first section. The second section lists results on representations of the braid groups $B_3$ and $B_4$ as well as results on the Temperley-Lieb algebras. This allows us to determine the eigenvalues of the braiding morphism $c_{V,V}$ in terms of a parameter $q^2$ in Section 3; here $V$ is an object in a tensor category of type $G_2$ or $G_{2,k}$ corresponding to the simple 7-dimensional representation of $G_2$. We determine the possible values of $q^2$ for a given tensor category of type $G_2$ or $G_{2,k}$ in the following section. Section 5 proves that $\End(V^{\otimes n})$ is already generated by the image of the braid group $B_n$ for all $n\in \N$, for categories of type $G_2$ and $G_{2,k}$. The already stated main result is then proved in Section 6. We also discuss connections of our approach  with the one proposed by Morrison, Peters and Snyder in \cite{MPS} using Kuperberg's spiders,
see \cite{Ku2} and subsequent work in \cite{MPS} and \cite{MST}, as well as possible generalizations to studying tensor categories of other exceptional types.

$Acknowledgement:$  H.W. would like to thank
UNC Wilmington and Lilit Martirosyan for supporting two visits, Eric Rowell for useful references and Pavel Etingof for advice concerning the symmetric case, which was missing in the first version of our paper. L.M. gratefully acknowledges a Reassignment Award from the College of Science and Engineering at UNC Wilmington. Both authors would like to thank MFO Oberwolfach for supporting us
as research fellows for two weeks in October 2024. Last, but not least, they would like to thank the referee for useful suggestions which improved the presentation of the material.

\section{Ribbon tensor categories of type $G_2$}

\subsection{Braided tensor categories}\label{braidedtensor}   
We recall a few basic facts about braided and ribbon
tensor categories, see e.g. \cite{EGNO}, \cite{Ks}, \cite{Turaev} for more details.
This serves mostly as motivation for the definitions in the next subsection.
In the following, we always deal with $\C$-linear tensor categories, i.e. all Hom spaces are complex vector spaces. We call a tensor category $\Ca$ semisimple if $\End(X)$ is a semisimple algebra for all objects $X$ in $\Ca$.

A braided tensor category $\Ca$ has canonical isomorphisms
$c_{V,W}: V\otimes W\to W\otimes V$ for any
object $V$, $W$ in $\Ca$. They satisfy the condition
\begin{equation}\label{braid1}
c_{U, V\otimes W}\ =\ (1_V\otimes c_{U,W}) (c_{U,V}\otimes 1_W),
\end{equation}
and a similar identity for $c_{U\otimes V, W}$. 
One can show that we obtain a representation of the braid group $B_n$
into $\End(V^{\otimes n})$ for any object $V$ in $\Ca$ via the map
\begin{equation}\label{braid2}
\sigma_i\ \mapsto 1_{i-1}\otimes  c_{V,V} \otimes 1_{n-1-i},
\end{equation}
where $1_k$ is the identity morphism on $V^{\otimes k}$. Recall that the 
 braid group $B_n$ is the group defined by generators $\sigma_i$, $1\leq i<n$ 
and relations $\sigma_i\sigma_{i+1}\sigma_i=
\sigma_{i+1}\sigma_i\sigma_{i+1}$ as well as $\sigma_i\sigma_j=\sigma_j
\sigma_i$ for $|i-j|\geq 2$. Let $\Delta_n\in B_n$ be defined inductively by 
$\Delta_2=\sigma_1$ and
$\Delta_n=\Delta_{n-1}\sigma_{n-1}\sigma_{n-2}\ ...\ \sigma_1$. 
Then it is well-known that $\Delta_n^2=(\sigma_1\sigma_2\ ...\ \sigma_{n-1})^n$
generates the center of $B_n$, see e.g. \cite{KT} Section 1.3.3. We shall also need the simple formula
\begin{equation}\label{braidflip}
\sigma_1\Delta_4=\Delta_4\sigma_3,\quad \sigma_3\Delta_4=\Delta_4\sigma_1.
\end{equation}
If $\Ca$ is a braided tensor category, an associated ribbon structure is given by maps $\Theta_W: W\to W$, $W\in$ Ob$(\Ca)$, which satisfy
\begin{equation}\label{ribbon}
\Theta_{V\otimes W}\ =\ c_{W,V}c_{V,W} (\Theta_V\otimes \Theta_W).
\end{equation}
Then we can prove by induction on $n$
that 
\begin{equation}\label{twistbraid}
\Theta_{V^{\otimes n}}= \Delta_n^2 \Theta_V^{\otimes n}.
\end{equation}
 If $V_\la$ is a simple object, 
the ribbon map just acts via a scalar, which we will denote by $\Theta_\la$.  
Let us also assume that the representation of $B_n$ into $\End(V^{\otimes n})$
is semisimple. Then the central element $\Delta_n^2$ acts in the simple
component labeled by $\la$ via a scalar denoted by $\al_{\la, n}$.
If the $B_n$-representation labeled by $\la$ acts nontrivially on the
$\End(V^{\otimes n})$-module $\Hom(V_\la,V^{\otimes n})$,  then it follows
from Eq \ref{twistbraid} that
\begin{equation}\label{scalarmatch}
\Theta_\la  =\ \Theta_V^n \alpha_{\la, n} ,
\end{equation}
where we identified $\Theta_V$ with the scalar via which it acts on $V$. We also assume that our tensor categories are semisimple and $rigid$, see \cite{Ks}, \cite{BK}. This, together with the ribbon structure, implies the existence of a canonical trace functional $Tr_X$ for $\End(X)$, for any object $X$ in our category. We will need  the following elementary result (see e.g. \cite{BK}, Lemma 2.4.1):

\begin{lemma}\label{nonzero trace} Let $X$ be a simple object in a rigid, semisimple ribbon tensor category. Then $\dim X:= Tr_X(1_X)\neq 0$.
\end{lemma}

\subsection{Categories of Type $G_2$}  

\begin{definition}\label{typedefinition} Let $G$ be a compact group.
We say a ribbon tensor category $\Ca$ is of type $G$ if it is a rigid semisimple
category for which there exists a 1-1 correspondence between 
isomorphism classes of simple objects of Rep $G$ and of $\Ca$ which defines an isomorphism
between their representation rings.
If $G$ is a compact Lie group of type $X_N$ for a Dynkin diagram $X_N$,
we may also just say that $\Ca$ is of type $X_N$.
\end{definition}
We will be particularly interested in ribbon tensor categories of type $G_2$.  
We first recall some basic facts about the roots and weights of the simple Lie algebra of type $G_2$.

With respect to the orthonormal unit vectors $\eps_1$, $\eps_2$, $\eps_3$ of $\R^3$, the roots of $\g$ can be written as $\Phi=\pm \{\eps_1-\eps_2,\eps_2-\eps_3, \eps_1-\eps_3, 2\eps_1-\eps_2-\eps_3, 2\eps_2-\eps_1-\eps_3, 2\eps_3-\eps_1-\eps_2\}$. The base can be chosen $\Pi=\{\al_1=\eps_1-\eps_2, \al_2=-\eps_1+2\eps_2-\eps_3\}$. The Weyl vector is given by $\rho=2\eps_1+\eps_2-3\eps_3$ and the Weyl group is $D_6$. The fundamental dominant weights are given by $\Delta=\{\La_1=\eps_1-\eps_3=(1,0,-1), \La_2=\eps_1+\eps_2-2\eps_3 = (1,1,-2)\}$,
and the set of dominant integral weights is given by
$$P_+=P(G_2)_+=\{ m_1\La_1+m_2\La_2, m_i\geq 0\}\ =\ \{ (\mu_1,\mu_2,-\mu_1-\mu_2), \mu_1
\geq\mu_2\geq 0\},$$
where $m_i$ and $\mu_i$ are non-negative integers for $i=1,2$. Our definitions imply  $\mu_1=m_1+m_2$ and $\mu_2=m_2$. 
It follows from Weyl's dimension formula that the $q$-dimension of the $U_q\g(G_2)$-module $V_\mu$ 
with highest weight  $\mu=(\mu_1,\mu_2, -\mu_1-\mu_2)$ 
is equal to
\begin{equation}\label{dimensionformula}
d_\mu\ =\ \frac{[\mu_1-\mu_2+1][2\mu_1+\mu_2+5][\mu_1+2\mu_2+4][3\mu_1+6][3\mu_2+3][3(\mu_1+\mu_2)+9]}
{[1][5][4][6][3][9]}.
\end{equation}
\begin{corollary} \label{dimestimate} Assume $q^2$ is not a primitive $\ell$-th root of unity, for $\ell\in \{ 2,3,4,5,6,9\}$. Then $d_\mu\neq 0$ for all $\mu$ with $|\mu|=\mu_1+\mu_2\leq n$ implies $[m]\neq 0$ for $1<m\leq 2n+5$ and $[3m]\neq 0$ for $1<m\leq n+3$.
\end{corollary}
$Proof.$ By our conditions on $q$, $d_\mu=0$ only if one of the factors in the numerator is equal to 0. The claim now follows from \ref{dimensionformula} by induction on $n$.

\subsection{Tensor product  rules for $G_2$}\label{tensorG2}   We will need to know how to decompose the tensor product of an irreducible
representation $V_\la$ with $V=V_{\Lambda_1}$, the simple representation of dimension 7.
We will review this here for the reader's convenience. Using these rules,
we can describe the sequence $\End(V^{\otimes i})\subset \End(V^{\otimes i+1})$ 
for $0\leq i<4$ by the following Bratteli diagram:

\begin{tikzpicture}[scale=.4]

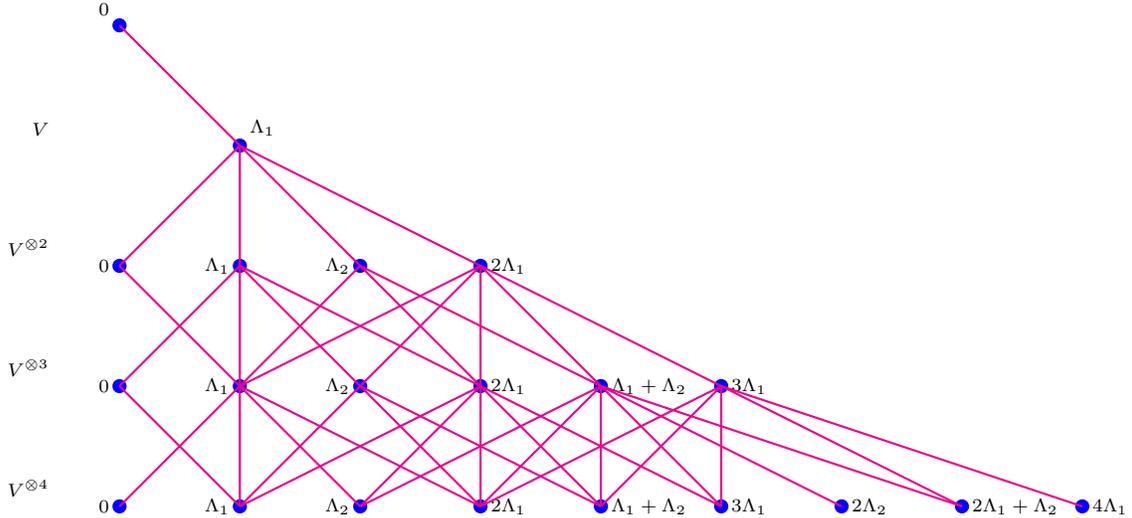
    \draw (0,-16) node[anchor=east]  {Bratteli diagram for $V^{\otimes n}$};
        \draw (-10, 4) node[anchor=south east]  {\tiny $0$};

    \draw [thick,blue,fill=blue](-10 cm,4 cm) circle (.2cm);

     \draw (-6, 0) node[anchor=south west]  {\tiny $\La_1$};
     \draw (-12, 0) node[anchor=south east]  {\tiny $V$};
\foreach \x in {2,...,4}      \draw (-12, -4*\x+4) node[anchor=south east]  {\tiny $V^{\otimes \x}$};

\foreach \x in {1,...,3}     \draw (-10, -4*\x) node[anchor= east]  {\tiny $0$};
\foreach \x in {1,...,3}     \draw (-6, -4*\x) node[anchor= east]  {\tiny $\La_1$};
\foreach \x in {1,...,3}     \draw (-2, -4*\x) node[anchor= east]  {\tiny $\La_2$};
\foreach \x in {1,...,3}     \draw (2, -4*\x) node[anchor= west]  {\tiny $2\La_1$};
\foreach \x in {2,3}     \draw (6, -4*\x) node[anchor= west]  {\tiny $\La_1+\La_2$};
\foreach \x in {2,3}     \draw (10, -4*\x) node[ anchor= west]  {\tiny $3\La_1$};
\draw (14, -12) node[anchor= west]  {\tiny $2\La_2$};
 \draw (18, -12) node[anchor= west]  {\tiny $2\La_1+\La_2$};
  \draw (22, -12) node[anchor= west]  {\tiny $4\La_1$};

    \draw [thick,blue,fill=blue](-6 cm,0 cm) circle (.2cm);
 \foreach \x in {0,...,3}
   \draw[xshift=\x cm,thick,blue,fill=blue] (3*\x cm-10 cm,-4 cm) circle (.2cm);

 \foreach \x in {0,...,5}
   \draw[xshift=\x cm,thick,blue,fill=blue] (3*\x cm-10 cm,-8 cm) circle (.2cm);
\foreach \x in {0,...,8}
\draw[xshift=\x cm,thick,blue,fill=blue, outer sep=0.3pt,
    inner sep=0.5pt] (3*\x cm-10 cm,-12 cm) circle (.2cm);

      \draw[thick, magenta] (-10,4) -- +(4, -4);

    	  \draw[thick, magenta] (-6,0) -- +(8, -8);
    	  \draw[thick, magenta] (-6,0) -- +(16, -8);
    	  \draw[thick, magenta] (-6,0) -- +(0, -12);
	  \draw[thick, magenta] (-10,-4) -- +(4, 4);
	  \draw[thick, magenta] (-6,-4) -- +(8, -4);
	  \draw[thick, magenta] (-6,-4) -- +(4, -4);

	  \draw[thick, magenta] (-6,-8) -- +(8, -4);

    	  \draw[thick, magenta] (-10,-4) -- +(8, -8);
    	  \draw[thick, magenta] (-10,-8) -- +(4, -4);
    	  \draw[thick, magenta] (-10,-8) -- +(4, 4);
    	  \draw[thick, magenta] (-10,-12) -- +(8, 8);
    	  \draw[thick, magenta] (-2,-8) -- +(4, -4);
    	  \draw[thick, magenta] (-2,-8) -- +(8, -4);
    	  \draw[thick, magenta] (-2,-4) -- +(8, -4);
    	  \draw[thick, magenta] (-2,-8) -- +(-4, -4);
    	  \draw[thick, magenta] (2,-4) -- +(0, -8);
    	  \draw[thick, magenta] (2,-4) -- +(-8, -4);
    	  \draw[thick, magenta] (2,-4) -- +(4, -4);
    	  \draw[thick, magenta] (2,-8) -- +(-8, -4);
    	  \draw[thick, magenta] (2,-8) -- +(8, -4);
    	  \draw[thick, magenta] (2,-8) -- +(4, -4);
    	  \draw[thick, magenta] (2,-8) -- +(-4, -4);
    	  \draw[thick, magenta] (2,-4) -- +(-4, -4);
    	  \draw[thick, magenta] (6,-8) -- +(-4, -4);
    	  \draw[thick, magenta] (6,-8) -- +(0, -4);
    	  \draw[thick, magenta] (6,-8) -- +(-8, -4);
    	  \draw[thick, magenta] (6,-8) -- +(4, -4);
    	  \draw[thick, magenta] (6,-8) -- +(8, -4);
    	  \draw[thick, magenta] (6,-8) -- +(12, -4);
    	  \draw[thick, magenta] (10,-8) -- +(0, -4);
    	  \draw[thick, magenta] (10,-8) -- +(12, -4);
    	  \draw[thick, magenta] (10,-8) -- +(-4, -4);
    	  \draw[thick, magenta] (10,-8) -- +(8, -4);
    	  \draw[thick, magenta] (10,-8) -- +(-8, -4);

\end{tikzpicture}

Let $\tilde V$ be a module of the semisimple Lie algebra $\g$ for which one explicitly knows the multiplicity $m_\om$ of each weight $\om$. If $V_\la$ is a simple $\g$ module with highest weight $\la$, it is a well-known consequence of Weyl's character formula that one can calculate the decomposition of the tensor product $V_\la\otimes \tilde V$ via the formula
%\begin{equation}\label{tensorrule} V_\la\otimes \tilde V\ \cong\ \bigoplus (\sum_{ w(\la+\om+\rho)=\mu+\rho}) {\rm sign}(w) m_\om)\ V_\mu, \end{equation}

\begin{equation}\label{tensorrule}  
V_\la\otimes\tilde V\ \cong\
\bigoplus_{\mu\in P_+}
\left(
\sum_{\substack{w\in W,\ \om\in\wt(\tilde V)\\
w(\la+\om+\rho)=\mu+\rho}}
\operatorname{sign}(w)m_\om
\right)V_\mu,
\end{equation}
where $w\in W$, the Weyl group of $\g$. 
Recall that for $\g=\g(G_2)$ the weights of $V=V_{\La_1}$ are
the short roots of $\g$ together with the zero weight. In this case the decomposition of the tensor product 
$V_\la\otimes V$ can
be described as follows (see e.g. \cite{MW}, Prop. 2.1 and Remark 2.2): 
Consider the hexagon centered at $\la+\rho$
and with corners $\la+\rho+\om$, with $\om$ running through
the short roots of $\g$. If this hexagon is contained in the dominant
Weyl chamber $C$, then $V_\la\otimes V$ decomposes into the direct sum
of irreducible $\g$-modules $V_\mu$ with $\mu+\rho$ running through the corners and the center of
the hexagon. If it is not contained in $C$, leave out all
the corners of the hexagon which are not in $C$; moreover, if
$\la=b\La_2$, also leave out $\la$ itself.

We can similarly visualize the weights of the adjoint representation $V_{\La_2}$ by
adding the long roots to the hexagon, and keeping track of the fact that the zero weight appears with multiplicity 2.
Using Eq. \ref{tensorrule}, or the same pictorial calculus as for the first fundamental representation and Young diagram notation for the weights
(i.e. $\La_1=[1,0]$ and $\La_2=[1,1]$), we obtain
$$[2,0]\otimes [1,1]\ =\  [3,1]\oplus [3,0]\oplus [2,1]\oplus [2,0]\oplus  [1,1]\oplus  [1,0].$$
Using the same procedure, we also obtain
$$[1,1]\otimes [1,1]\ =\  [2,2]\oplus [3,0]\oplus [2,0]\oplus  [1,1]\oplus  [0,0].$$
Using these results, the tensor product rules for $[2,0]\otimes [1,0]^{\otimes 2}$ and (formally) $[2,0]=[1,0]^{\otimes 2}-[1,1]-[1,0]-[0,0]$, we calculate
$$[2,0]\otimes [2,0]\ =\  [4,0]\oplus [3,1]\oplus [2,2]\oplus [3,0]\oplus 2[2,1]\oplus 2[2,0]\oplus  [1,1]\oplus  [1,0]\oplus [0,0].$$
We finally also observe that $V_\la\otimes V_{\La_2}$ with $\la=[a+1,a]$ does not contain the summand $V_{[a+1,a+1]}$ by Eq. \ref{tensorrule}. Here the weight $\la+\rho+[0,1]=[a+3,a+2]$ cancels with $w(\la+\rho+[-1,2])$ in Eq. \ref{tensorrule}, where $w$ permutes the coordinates and $\rho=[2,1]$. One similarly shows that $V_{[a+1,a]}\not\subset V_{[a,a]}\otimes V_{\La_2}$.

Let $|\la|=\la_1+\la_2$ for any element in the weight lattice of $G_2$. Then it is easy to see that
a highest weight module $V_\la$ appears in the $|\la|$-th tensor power of $V$ for the first time. We also have
the following consequence of the tensor product rules
\begin{equation}\label{tensoradjoint}
V_{\la}\otimes V_{\La_2}\ =\ V_{\la+\La_2}\oplus\ \bigoplus V_\mu,
\end{equation}
where $|\mu|\leq |\la|+1$. Indeed, this follows from the well-known general fact that if $V_\mu\subset V_{\la}\otimes V_{\La_2}$, then $\mu=\la+\om$ for some weight $\om$ in $V_{\La_2}$ (see e.g. \cite{MW}, Prop. 2.1 and Remark 2.2 for more details). As $V_{\La_2}$ is the adjoint representation, its weights are the roots and the zero weight. So, the claim follows from the fact that $|\om|<2$ for all weights $\om\neq \La_2$ of $V_{\La_2}$.

The action of $B_n$ on $V^{\otimes n}$ also gives us representations of $B_n$ on $\Hom(V_\mu, V^{\otimes n})$ for all $V_\mu\subset V^{\otimes n}$.  Then our tensor product rules imply
the restriction rule 
\begin{equation}\label{braidrestrict}
\Hom(V_\mu, V^{\otimes n})_{|B_{n-1}}\ \cong\ \oplus_{\nu}\Hom(V_\nu,V^{\otimes n-1}),
\end{equation}
where the summation on the right goes over all $\nu$ such that $V_\nu\subset V^{\otimes n-1}$ and $V_\mu\subset V_\nu\otimes V$.

\subsection{Tensor  categories of type $G_{2,k}$}\label{fusion:sec} Let \(\ell=k+12\). For each integer \(k\geq -2\) and each primitive \(\ell\)-th root of unity \(q^2\), we associate a fusion tensor category \(\bar{\U}_q\). It can be obtained as a quotient of the category of tilting modules of the Drinfeld-Jimbo quantum group $\U_q$ for $q^2$ a primitive $(k+12)$-th root of unity, see e.g. \cite{AP}, \cite{Sawin}.

\begin{center}
  \begin{tikzpicture}[scale=.4]

    % Title with ell = 13
    \draw (0,-9) node[anchor=east] {Weyl alcove for $\mathfrak{g}$, $k=1$, $\ell = 13$};

    % Full lattice points (gray)
    \foreach \y in {4,...,6}
      \foreach \x in {0,...,\y}
        \draw[xshift=\x cm, thick, gray] (\x cm, 3*\x cm) circle (.2cm);

    % Blue points (levels 0 to 3)
    \foreach \y in {0,...,3}
      \foreach \x in {0,...,\y}
        \draw[xshift=\x cm, thick, blue, fill=blue] (\x cm, 3*\x cm) circle (.2cm);

    % Level 4
    \foreach \x in {0,...,2}
      \draw[xshift=\x cm, thick, blue, fill=blue] (\x cm, 3*\x cm + 6 cm) circle (.2cm);
    \foreach \x in {3,...,4}
      \draw[xshift=\x cm, thick, gray] (\x cm, 3*\x cm + 6 cm) circle (.2cm);

    % Level 5
    \draw[xshift=0cm, thick, blue, fill=blue] (0cm, 3*0cm + 12 cm) circle (.2cm);
    \foreach \x in {1,...,2}
      \draw[xshift=\x cm, thick, gray] (\x cm, 3*\x cm + 12 cm) circle (.2cm);

    % Level 6
    \foreach \x in {0,...,0}
      \draw[xshift=\x cm, thick, gray] (\x cm, 3*\x cm + 18 cm) circle (.2cm);

    % Weight labels
    \foreach \x in {2,...,6}
      \draw[xshift=\x cm] (\x cm, 3*\x cm) node [anchor=north west] {\tiny $\x\Lambda_1$};
    \foreach \x in {2,...,4}
      \draw[xshift=\x cm] (\x cm, 3*\x cm + 6 cm) node [anchor=north] {\tiny $\x\Lambda_1 + \Lambda_2$};
    \foreach \x in {2,...,2}
      \draw[xshift=\x cm] (\x cm, 3*\x cm + 12 cm) node [anchor=north] {\tiny $\x\Lambda_1 + 2\Lambda_2$};
    \foreach \y in {2,...,3}
      \draw[yshift=\y cm] (0cm, 5*\y cm) node [anchor=north east] {\tiny $\y\Lambda_2$};
    \foreach \y in {2,...,2}
      \draw[yshift=\y cm] (2.5 cm, 5*\y cm + 3cm) node [anchor=north] {\tiny $\Lambda_1 + \y\Lambda_2$};

    % Origin and fundamental weights
    \draw (0 cm, 0 cm) node[anchor=north] {\tiny $0$};
    \draw (0 cm, 6 cm) node[anchor=north east] {\tiny $\Lambda_2$};
    \draw (2.2 cm, 3 cm) node[anchor=north west] {\tiny $\Lambda_1$};
    \draw (2.2 cm, 9 cm) node[anchor=north] {\tiny $\Lambda_1 + \Lambda_2$};

    % Original chamber boundaries (dotted)
%    \draw[dotted, thick] (0,0) -- +(14, 21);
%    \draw[dotted, thick] (0,0) -- +(0, 21);

    % Decorative elements
%    \draw[dotted, thick, magenta] (2,3) -- +(0, 18);
 %   \draw[dotted, thick, magenta] (2.1,3.1) -- +(12, 18);
    \draw[thick, red] (-2,17) -- +(12.8, -6.4);

% Define custom loop style with stealth arrow
\tikzset{
  looparrow/.style={
    thick,
    orange,
    postaction={
      decorate,
      decoration={
        markings,
        mark=at position 0.6 with {
          \arrow[scale=0.55, orange]{Stealth}
        }
      }
    }
  }
}

%\draw[thick, orange, ->, >=Stealth]
 % (2,9) .. controls +(-1.2,0.8) and +(-1.2,-0.8) .. (2,9);

\draw[thick, orange, ->, >=Stealth]
  (6,9) .. controls +(1.8,1.2) and +(-0,1.2)  .. (6,9);

\draw[thick, orange, ->, >=Stealth]
  (2,3) .. controls +(-0.8,-1.2) and +(0.4,-1.2) .. (2,3);

\draw[thick, orange, ->, >=Stealth]
  (4,6) .. controls +(-0.8,-1.2) and +(0.4,-1.2) .. (4,6);
  
\draw[thick, orange, ->, >=Stealth]
  (2,9) .. controls +(-1.2,1.2) and +(-1.2,-0.3) .. (2,9);

   % Boundary lines
    \draw[thick, black] (-2,-7.5) -- (-2, 17);         % vertical line
   \draw[dotted, thick, black] (-2,17) -- (-2, 21);         % vertical line
    \draw[dotted,thick, black] (-2,-7.5) -- (18, 21);          % diagonal   
   
\draw[thick, black] (-2,-7.5) -- (10.74,10.65);

\draw[thick, orange] (0,6) -- (4,6);
\draw[thick, orange] (2,9) -- (6,9);
\draw[thick, orange] (4,6) -- (6,9);
\draw[thick, orange] (0,6) -- (2,9);
\draw[thick, orange] (0,0) -- (4,6);
\draw[thick, orange] (0,12) -- (4,6);
\draw[thick, orange] (0,12) -- (4,12);
\draw[thick, orange] (0,6) -- (2,3);
\draw[thick, orange] (6,9) -- (4,12);
\draw[thick, orange] (2,9) -- (4,12);
\node[gray, thick, rotate=90] at (0,20) {$\cdots$};
\node[gray, thick, rotate=55] at (1.5,20) {$\cdots$};
\node[gray, thick, rotate=55] at (5.5,20) {$\cdots$};
\node[gray, thick, rotate=55] at (9.5,20) {$\cdots$};
\node[gray, thick, rotate=55] at (13.5,20) {$\cdots$};
  \end{tikzpicture}
\end{center} 
Recall that $\La_2$ and $\ep_1-\ep_3$ are the highest long and short roots. The simple objects of $\bar{\U}_q$ are labeled by the sets
$P_{+,k}\subset P_+$ by
$$P_{+,k}\ =\ \begin{cases}\ \{ \la\in P_+,\ (\la+\rho,\La_2)=3(\la_1+\la_2+3)< k+12\}, &{\rm if}\ 3|k,\cr 
\ \{ \la\in P_+,\ (\la+\rho,\ep_1-\ep_3)=2\la_1+\la_2+5< k+12\}, &{\rm if}\ 3\nmid k.
\end{cases}
$$
{\bf Warning:} Our number $k$ is three times what is usually called the level for integrable representations of affine Kac-Moody algebras of type $G_2$.  
Here are two examples of $P_{+,k}$.
 The following graph describes $P_{+,k}$ for $k=9$, or, equivalently, $\ell=21$.

\begin{center}
  \begin{tikzpicture}[scale=.4]

    % Title with ell = 21
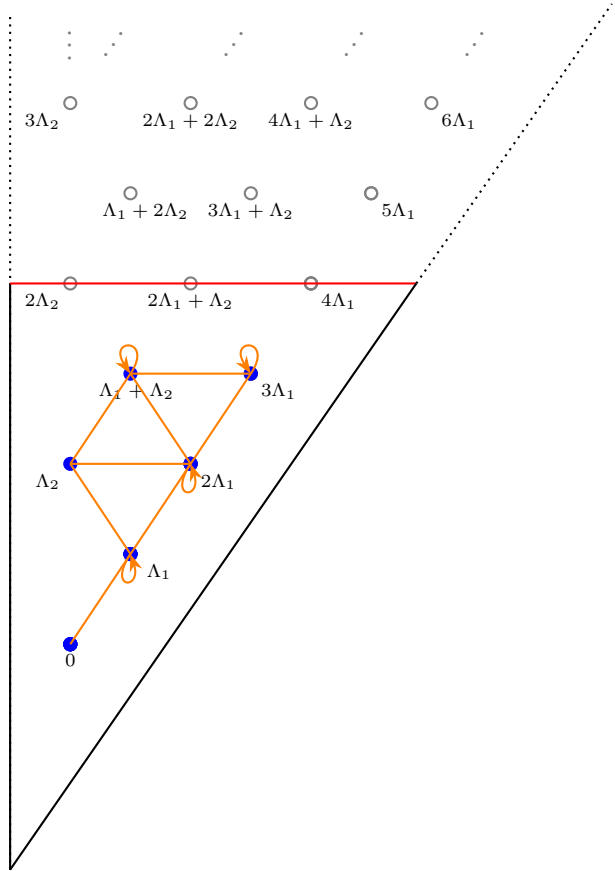
    \draw (0,-9) node[anchor=east] {Weyl alcove for $\mathfrak{g}$, $\ell = 21$};

    % Full lattice points (gray)
    \foreach \y in {4,...,6}
      \foreach \x in {0,...,\y}
        \draw[xshift=\x cm, thick, gray] (\x cm, 3*\x cm) circle (.2cm);

    % Blue points (levels 0 to 3)
    \foreach \y in {0,...,3}
      \foreach \x in {0,...,\y}
        \draw[xshift=\x cm, thick, blue, fill=blue] (\x cm, 3*\x cm) circle (.2cm);

    % Level 4
    \foreach \x in {0,...,1}
      \draw[xshift=\x cm, thick, blue, fill=blue] (\x cm, 3*\x cm + 6 cm) circle (.2cm);
    \foreach \x in {2,...,4}
      \draw[xshift=\x cm, thick, gray] (\x cm, 3*\x cm + 6 cm) circle (.2cm);

    % Level 5
    \draw[xshift=0cm, thick, gray] (0cm, 3*0cm + 12 cm) circle (.2cm);
    \foreach \x in {1,...,2}
      \draw[xshift=\x cm, thick, gray] (\x cm, 3*\x cm + 12 cm) circle (.2cm);

    % Level 6
    \foreach \x in {0,...,0}
     \draw[xshift=\x cm, thick, gray] (\x cm, 3*\x cm + 18 cm) circle (.2cm);

    % Weight labels
    \foreach \x in {2,...,6}
      \draw[xshift=\x cm] (\x cm, 3*\x cm) node [anchor=north west] {\tiny $\x\Lambda_1$};
    \foreach \x in {2,...,4}
      \draw[xshift=\x cm] (\x cm, 3*\x cm + 6 cm) node [anchor=north] {\tiny $\x\Lambda_1 + \Lambda_2$};
    \foreach \x in {2,...,2}
      \draw[xshift=\x cm] (\x cm, 3*\x cm + 12 cm) node [anchor=north] {\tiny $\x\Lambda_1 + 2\Lambda_2$};
    \foreach \y in {2,...,3}
      \draw[yshift=\y cm] (0cm, 5*\y cm) node [anchor=north east] {\tiny $\y\Lambda_2$};
    \foreach \y in {2,...,2}
      \draw[yshift=\y cm] (2.5 cm, 5*\y cm + 3cm) node [anchor=north] {\tiny $\Lambda_1 + \y\Lambda_2$};

    % Origin and fundamental weights
    \draw (0 cm, 0 cm) node[anchor=north] {\tiny $0$};
    \draw (0 cm, 6 cm) node[anchor=north east] {\tiny $\Lambda_2$};
    \draw (2.2 cm, 3 cm) node[anchor=north west] {\tiny $\Lambda_1$};
    \draw (2.2 cm, 9 cm) node[anchor=north] {\tiny $\Lambda_1 + \Lambda_2$};

    % Orange full circular arrows based at specified blue dots
 %   \draw[orange, thick, ->] (2,9) arc (0:359:0.5);
 %   \draw[orange, thick, ->] (6,9) arc (0:359:0.5);
 %   \draw[orange, thick, ->] (2,3) arc (0:359:0.5);
  %  \draw[orange, thick, ->] (4,6) arc (0:359:0.5);

    % === Boundary lines ===
 \draw[thick, black] (-2,-7.5) -- (-2,12);       % vertical stays at x = -2
  \draw[dotted, thick, black] (-2,-7.5) -- (-2,21);       % vertical stays at x = -2
\draw[thick, red] (-2,12) -- (11.5,12);         % red line shifted right by 2 cm
\draw[thick, black] (11.5,12) -- (-2,-7.5);     % diagonal shifted right, intersects vertical at same point
           % solid diagonal
\draw[dotted, thick, black] (11.5,12) -- (18, 21.3);           % continuation, dotted

  %  \draw[thick, red] (-2,19) -- +(12.8, -6.4);

  \draw[thick, orange] (0,6) -- (4,6);
\draw[thick, orange] (2,9) -- (6,9);
\draw[thick, orange] (4,6) -- (6,9);
\draw[thick, orange] (0,6) -- (2,9);
\draw[thick, orange] (0,0) -- (4,6);
%\draw[thick, orange] (0,12) -- (4,6);
%\draw[thick, orange] (0,12) -- (4,12);
\draw[thick, orange] (0,6) -- (2,3);
%\draw[thick, orange] (6,9) -- (4,12);
%\draw[thick, orange] (2,9) -- (4,12);
\draw[thick, orange] (4,6) -- (2,9);

% Define custom loop style with stealth arrow
\tikzset{
  looparrow/.style={
    thick,
    orange,
    postaction={
      decorate,
      decoration={
        markings,
        mark=at position 0.6 with {
          \arrow[scale=0.55, orange]{Stealth}
        }
      }
    }
  }
}

\draw[thick, orange, ->, >=Stealth]
  (2,9) .. controls +(0.8,1.2) and +(-0.8,1.2) .. (2,9);

\draw[thick, orange, ->, >=Stealth]
  (6,9) .. controls +(0.8,1.2) and +(-0.8,1.2) .. (6,9);

\draw[thick, orange, ->, >=Stealth]
  (2,3) .. controls +(-0.8,-1.2) and +(0.4,-1.2) .. (2,3);

\draw[thick, orange, ->, >=Stealth]
  (4,6) .. controls +(-0.8,-1.2) and +(0.4,-1.2) .. (4,6);

\node[gray, thick, rotate=90] at (0,20) {$\cdots$};
\node[gray, thick, rotate=55] at (1.5,20) {$\cdots$};
\node[gray, thick, rotate=55] at (5.5,20) {$\cdots$};
\node[gray, thick, rotate=55] at (9.5,20) {$\cdots$};
\node[gray, thick, rotate=55] at (13.5,20) {$\cdots$};
  \end{tikzpicture}
\end{center}
The decomposition of tensor products of simple objects in $\bar{\U}_q$ is given by the Kac-Walton formula, which can be described as follows. Let $W_k$ be the affine Weyl group generated by the simple reflections of the Weyl group $W$ of type $G_2$ and the additional reflection $s_0$ in the hyperplane given by $x_1+x_2= 4+k/3$ if $3|k$, and by $2x_1+x_2= k+12$ if $3\nmid k$.
Here we identify $\h^*$ with $(x_1,x_2,-x_1-x_2)\subset \R^3$.
The action is then given by 
$$s_0.\la= s_0(\la+\rho)-\rho = \la -\frac{2((\la+\rho, \theta)-(k+12))}{(\theta,\theta)}\theta, $$
where $\theta=\La_2=(1,1,-2)$ for $3|k$, and $\theta=\La_1=(1,0,-1)$ for $3\nmid k$.
The modified tensor product is defined by making the following adjustments to the usual tensor product for $G_2$ in \ref{tensorrule} (see e.g. \cite{Sawin} for more details).

(a) If $V_\mu$ appears in $V_\la\otimes V$ such that $w.\mu=\mu$ for some $1\neq w\in W_k$, we remove it.

(b) If $V_\mu$ appears in $V_\la\otimes V$ with $\mu\not\in P_{+,k}$,  such that $w.\mu\neq\mu$ for all $1\neq w\in W_k$,
then we replace $V_\mu$ by $\operatorname{sign}(w)V_{w.\mu}$ for the unique $w\in W_k$ such that $w.\mu\in P_{+,k}$.

The tensor product of $V_\la$ with $V$ in $\Ca(G_{2,k})$ can then be described
as follows:
If $3\nmid k$, we remove all those $\mu$ in \ref{tensorrule}  for which $2\mu_1+\mu_2\geq k+7$, and $V_\la$, if $2\la_1+\la_2=k+6$. So the tensor product of $V_\la$ with $V$ for $k=1$ is isomorphic to the direct sum $\oplus_\mu V_\mu$ with $\mu$ ranging over all vertices $\mu$ connected with $\la$ by a single edge or a loop in the picture above.  The graph on the previous page describes this for $k=1$, or, equivalently, $\ell=13$.

If $3|k$, we simply remove all the $V_\mu$ in \ref{tensorrule} for which $3(\mu_1+\mu_2)=k+3$.

\begin{definition}\label{fusiondef} We say a semisimple rigid ribbon tensor category $\Ca$ is of type $G_{2,k}$ if its simple objects are labeled by the elements of $P_{+,k}$ and it has the same fusion rules as $\bar{\U}_q$ for $q^2$ a primitive $(k+12)$-th root of unity.
\end{definition}

\subsection{Elementary properties} Recall that the Fibonacci category has two simple objects, 1 and $X$, with the only nontrivial fusion rule given by $X\otimes X=1\oplus X$. It was already shown in \cite{FK} (see also \cite{KW}) that there are only two categories with such fusion rules, up to equivalence. They correspond to sets $\{ \xi, \xi^{-1}\}$ of primitive fifth roots of unity, up to equivalence of tensor categories. 

Another well-known tensor category, $SO(3)_7$, has four simple objects labeled by $(1)$, $(3)$, $(5)$ and $(7)$. Its tensor product rules are given by the ones of $SO(3)$, by setting $(9)$ equal to the empty object, and reflecting objects $(11)$ to -$(7)$, and $(13)$ to -$(5)$; here $(i)$ labels the $i$-dimensional irreducible representation of SO(3). Then we have, e.g.
$$(5)\otimes (7)\ =\ (3)\oplus (5)\oplus (7)\oplus (9) \oplus (11)\ \cong\ (3)\oplus (5).$$
It follows from \cite{KW} and \cite{TW2} that tensor categories of type $SO(3)_{m-2}$ for $m$ odd
are classified up to equivalence of tensor categories by sets $\{ \xi, \xi^{-1}\}$ of primitive $m$-th roots of unity.
We would like to thank Eric Rowell for telling us about the connection between $G_{2,6}$ and $SO(3)_7$.

\begin{lemma}\label{lowtensor} (a) Let $\Ca$ be a tensor category of type $G_{2,\infty}=G_2$ or type $G_{2,k}$, with $k\geq 1$. We have $\dim \Hom(V_{\La_1+\La_2},V^{\otimes 4})=8$  except for $k\in\{ 3,6\}$. 

(b) Categories $\Ca$ of type $G_{2,3}$ and $G_{2,6}$ have the fusion rules of the Fibonacci category and of the category $SO(3)_7$, and categories of type $G_{2,0}$ are equivalent to $Vec$.

(c) (\cite{Ru}, Section 3.2) Categories $\Ca$ of type $G_{2,-1}$ have the fusion rules of $SO(3)_{11}$, and categories of type $G_{2,-2}$ have the fusion rules of the Deligne product of two  Fibonacci categories.

(d) Tensor categories of type $G_{2,k}$ with $k\in \{ -2,-1, 0, 3, 6\}$ are classified.
\end{lemma}

$Proof.$ Parts (a) and (b) are well known and can be checked explicitly from the given fusion rules. 
Part (c) was shown in \cite{Ru}, Section 3.2. Categories with the fusion rules of $SO(3)_{m-2}$ are classified in \cite{TW} and \cite{Cp}, which includes the Fibonacci category for $m=5$.

\subsection{Quantum Casimir}\label{quantumCasimir:sec} We denote by $C_\la$ the scalar via which the Casimir element acts on the irreducible representation $V_\la$. Using notations and normalizations as in Section \ref{typedefinition}, we have
$$C_\la=(\la+2\rho,\la).$$
Here \((-, -)\) denotes the normalized Killing form. In particular, \((\la+2\rho,\la)\in \mathbb Z\) for dominant weights \(\la\).
We shall need the following results about twists and quantum Casimirs, see
 \cite{Dr}, \cite{Turaev}
\begin{proposition}
\label{prop:eigenvalue}
Let $V_\la, V_\mu, V_{\La_1}=V$ be
simple $\U_q$-modules with highest weights
$\la, \mu, \La_1$ respectively, and such that $V_\mu$ is a submodule of
$V_\la\otimes V$. Let us write $c_{\la,\mu}$ for the braiding morphism
$V_\la\otimes V_\mu\to V_\mu\otimes V_\la$. Then,
$$(c_{{\La_1},{\la}}c_{{\la},{\La_1}})_{|V_\mu}=
q^{C_\mu- C_\la-C_V}1_{V_\mu},$$

where  $C_V=C_{\La_1}$.
Moreover, the twisting factors $\Theta_\la$ and the scalar $x_{\la,n}$ via which $\Delta_n^2$ acts on $V_\la\subset V^{\otimes n}$ are given by
$\Theta_\la\ =\ q^{C_\la}$ and by $x_{\la,n}=\Theta_\la\Theta_V^{-n}=q^{C_\la-nC_V}$.
\end{proposition}

\subsection{Path representations}\label{pathrep} We review a few well-known facts about path bases, see e.g. \cite{MW},\cite{MW2} for more details. For simplicity, we just describe this for the algebras $A_n=\End(V^{\otimes n})$ for $V$ our standard generating object for a tensor category of type $G_2$ or $G_{2,k}$. A path $t$ of length $n$ is a sequence of elements $\mu^{(i)}$, $0\leq i\leq n$ in $P_+$ or $P_{+,k}$
$$t:\quad \mu^{(0)}=[0,0]\  \to\  \mu^{(1)}\  \to \mu^{(2)}\  \to \ ...\ \to\  \mu^{(n)}$$
such that $V_{\mu^{(i+1)}}\subset V_{\mu^{(i)}}\otimes V$. We then obtain a well-defined path idempotent
$p_t=\prod_i z_{\mu^{(i)}}$, where $z_{\mu^{(i)}}$ is the central idempotent in $\End(V^{\otimes i})$
for the isotypical component for the module $V_{\mu^{(i)}}\subset V^{\otimes i}$.
One observes that the $B_n$ module $\Hom(V_\mu, V^{\otimes n})$ has a basis $(\phi_t)$ labeled by all paths $t$ of length $n$ such that $t{(n)}=\mu$, where $p_t \phi_t=\phi_t\neq 0$. If $S_i$ is the matrix representing the $i$-th braid generator $\sigma_i$, we have
$$S_i \phi_t=\sum_s a_{st}^{(i)}\phi_s,$$
where the summation goes over all paths $s$ which coincide with $t$ except possibly for the $i$-th element. As $S_i$ commutes with the image of $\End(V^{\otimes i-1})\otimes 1_2\otimes \End(V^{\otimes n-i-1})$, we can assume that the matrix coefficients of $S_i$ only depend on the diagrams $\mu^{(j)}$, $i-1\leq j\leq i+1$ of the given paths.
In particular, we can write $S_i$ as a direct sum of matrix blocks  $S_i(\mu^{(i-1)}, \mu^{(i+1)})$, whose rows and columns are labeled by the set $\P(\mu^{(i-1)}, \mu^{(i+1)})$ of paths of length 2 from $\mu^{(i-1)}$ to $\mu^{(i+1)}$. 
The same applies for each eigenprojection of  $S_i$, i.e. the projection onto an eigenspace corresponding to a fixed eigenvalue. 

\subsection{Affine braid group $AB_2$}\label{affine:sec} We will use the affine braid group $AB_2$ to calculate the matrix blocks $S_i=S_i(\la,\nu)$ of the previous subsection in certain cases. It is defined via generators $\tau$ and $\sigma$ and relation
$\sigma\tau\sigma\tau=\tau\sigma\tau\sigma$. First observe that we obtain a homomorphism $AB_2\to B_{n+1}$ via the maps
\begin{equation}\label{affinebraidrel}
\tau\mapsto \Delta_n^2,\hskip 3em \sigma\mapsto \sigma_n.
\end{equation}
This can be seen easier if one maps $\tau$ to $\Delta_n^2\Delta_{n-1}^{-2}=\sigma_{n-1}\sigma_{n-2}\ ...\ \sigma_1^2\sigma_2\ ...\ \sigma_{n-1}$. But as $\Delta_{n-1}^2$ commutes with both $\sigma_n$ and $\Delta_n^2$, the braid relation  also holds for the assignments in Eq.\ref{affinebraidrel}.

We recall how one can obtain representations of $AB_2$ from representations of ordinary braid groups using
(quantum) Jucys-Murphy elements. This was first done for Hecke algebras in \cite{Cherednik}.

\begin{lemma}\label{affinebraidlemma} Let $\la$ and $\nu$ be labels in our $G_2$ Bratteli diagram at level $n-1$ and $n+1$. Then we obtain a representation of $AB_2$ on the vector space $W(\la,\nu)$ with a basis labeled by the set $\P(\la,\nu)$ of paths $t:\la\to\mu_t\to\nu$ of length 2 given by
$$\tau\mapsto T: t\mapsto x_t\ t\hskip 3em \sigma\mapsto A=S_n(\la,\nu),$$
with $x_t=q^{C_{\mu_t}-nC_V}$ and with notations as in Sections \ref{quantumCasimir:sec} and \ref{pathrep}. Then $ATAT$ acts via the scalar $\delta=q^{C_\la+C_\nu-2nC_V}$ on  $W(\la,\nu)$. Moreover, 
if \(a_{rs}\) denotes the \((r,s)\)-entry of the matrix \(A\), then
\begin{equation}\label{Ainverse}
(A^{-1})_{rs}=\delta^{-1} x_rx_s a_{rs}.
\end{equation}
\end{lemma}
%q^{C_\la+C_\nu-C_{\mu_r}-C_{\mu_s}}

$Proof.$ The first statement just follows from Eq \ref{affinebraidrel} and the definitions of the path representations. Using the remark below Eq \ref{affinebraidrel}, we see that $\sigma_n\Delta_n^2\Delta_{n-1}^{-2}\sigma_n=\Delta_{n+1}^2\Delta_n^{-2}$, from which one easily deduces that $TATA=\Delta_{n+1}^2\Delta_{n-1}^2$. One obtains Eq \ref{Ainverse} by multiplying this equation by $A^{-1}$.

\medskip

We consider representations of $AB_2$ on a vector space $W$ with a basis of eigenvectors $\{v_r\}$ with {\it mutually distinct} eigenvalues $x_r$ of $T$ such that $ATAT$ acts via a scalar $\delta$. The following lemma is a version of a well-known argument using quantum Jucys-Murphy elements. As we did not find a reference of it for our given version, we provide a proof for the reader's convenience.

\begin{lemma}\label{2by2matrixlem} Assume that $(A-\la_1)(A-\la_2)=0$. 

(a) We can have a 1-dimensional subrepresentation of $AB_2$ only if $\delta=x_r^2\la_i^2$ for suitable indices $r$ and $i\in \{ 1,2\}$.  

(b) Any indecomposable subrepresentation of $W$ has at most dimension 2, with the basis $\{ v_r, v_s\}$ for some indices $1\leq r,s\leq n$. It is irreducible if $\la_1/\la_2\neq -(x_r/x_s)^{\pm 1}$. \end{lemma}

$Proof.$ Part (a) follows from the definitions. Since $(A-\la_1)(A-\la_2)=0$, we have $ A+\la_1\la_2A^{-1}=(\la_1+\la_2)I$.  Together with Eq.~\ref{Ainverse}, this gives, for $r\neq s$, $ (\delta+\la_1\la_2x_rx_s)a_{rs}=0. $ Thus, if $a_{rs}\neq 0$, then $ x_rx_s=-\frac{\delta}{\la_1\la_2}. $ Since the $x_r$ are mutually distinct, for each $s$ there is at most one $r\neq s$ for which this can occur. Thus $W$ is a direct sum of subrepresentations of dimension at most two. Using e.g. the proof of \cite{MW}, Lemma 1.8, we calculate the diagonal entries of the element $(A-\la_11)$ as 
\begin{equation}\label{2dimformula} (A-\la_11)_{rr}=- \frac{\la_2x_s+\la_1x_r}{(x_r-x_s)}, \hskip 3em (A-\la_11)_{ss}=\frac{\la_2x_r+\la_1x_s}{(x_r-x_s)}. 
\end{equation} 
We see that they are both nonzero under the conditions stated in (b). As $A-\la_1$ has at most rank 1 on this subspace, the same applies to its off-diagonal elements. 

\subsection{Irreducibility criteria} We define an equivalence relation among paths as the transitive closure of $s\leftrightarrow t$ if 
$a^{(i)}_{st}a^{(i)}_{ts}\neq 0$ for the matrix entries of either a braid matrix $S_i$, or of one of its eigenprojections. It is well-known, and easy to show that a module $V_\mu$ is irreducible if all its paths are equivalent. Here is an easy example which will be used later.
\begin{lemma}\label{path:lem}
Assume the $B_{n-1}$ module $\Hom(V_\la, V^{\otimes n-1})$ and the $B_n$ modules $\Hom(V_\mu, V^{\otimes n})$
are simple for all $\mu$ such that $V_\mu$ is a summand of $V_\la\otimes V$, and also assume that the dimensions of $V_\la$ and all the involved $V_\mu$ are nonzero. Then also the $B_{n+1}$ module $\Hom(V_\la, V^{\otimes n+1})$ is simple.
\end{lemma}

 $Proof.$ Let $P$ be the projection onto $\1\subset V^{\otimes 2}$, and let $P_n$ be the corresponding eigenprojection of $S_n$. Then it is well-known (see e.g. \cite{Wexc}, Proposition 4.10) that its diagonal entries are given by $d_\mu/d_\la d_V$
 for the paths of length $n+1$ ending in $  ... \to \la\to\mu\to\la$. As $P_n$ has rank 1 for a given block, also its off-diagonal entries are nonzero. By assumption, the $B_n$ submodules $\Hom(V_\mu, V^{\otimes n})$ are irreducible. The nonzero off-diagonal entries of $P_n$ provide maps between these submodules.

 \begin{remark}\label{newpath} (a) Let $V_\ga\subset V^{\otimes 2}$ and $\la_\ga$ the eigenvalue of $c_{VV}$ by which it acts on $V_\ga$. Then it follows from our definitions that the multiplicity of $\la_\ga$ of the matrix block $S_i(\la,\nu)$ (i.e. we have $\mu^{(i-1)}=\la$ and $\mu^{(i+1)}=\nu$) is equal to the multiplicity of $V_\nu$ in $V_\la\otimes V_\ga$.
 
(b)  A similar argument as in the proof of Lemma \ref{path:lem} will also be used in the following context (see e.g. \cite{MW2}, Section 2.4 for more details). Assume  $\la=\mu^{(i-1)}$ and  $\nu=\mu^{(i+1)}$  are weights in a path such that $|\nu|=|\la|+1$ and $\la\subset\nu$. By Eq. \ref{tensorrule} and (a), it follows that the eigenprojection $P_i$ of $S_i$ corresponding to $\la_{\La_1}$ has rank 1.
 Assume that $P_i$ has nonzero diagonal entries with respect to our path basis. If $\Hom(V_\mu, V^{\otimes i})$ are irreducible $B_i$ modules for all weights $\mu$ appearing in a path $\la\to \mu\to \nu$, then all these submodules must belong to the same $B_{i+1}$ submodule of $\Hom(V_\nu, V^{\otimes i+1})$. %This argument was already used in \cite{MW} to show that $B_n$ generates $\End(V^{\otimes n})$ for tensor categories of type $G_2$, under some additional conditions. It will again be used in this paper in a more general setting.
 \end{remark}
\section{Representations of $B_3$ and $B_4$}

\subsection{Definitions}\label{definitionsection}
 We denote by $K_n$ the quotient of the group algebra $\C B_n$ of the braid group $B_n$
subject to the relations $(\sigma_i-\la_1)(\sigma_i-\la_2)(\sigma_i-\la_3)=0$,
for all $i$, where the parameters $\la_1,\la_2,\la_3\in\C$ are fixed. It was shown by Coxeter \cite{Cox}
that $B_n$ with the additional relation $\sigma_i^3=1$ is a finite group only if $n\leq 5$. The corresponding finite quotient groups are exceptional complex reflection groups labeled by $G_4$, $G_{25}$ and $G_{32}$ in \cite{ST}. The following
theorem is a nontrivial generalization of it.

\begin{theorem} \cite{MHecke} The algebras $K_n$, $n=3,4,5$, are flat deformations of the group algebras of the complex reflection groups $G_4$, $G_{25}$ and $G_{32}$.  In particular, they have bases with respect to which the standard generators
act via matrices whose entries are Laurent polynomials in the eigenvalues $\la_i$ over a finite extension of $\Z$. 
\end{theorem}
The structure of the algebra $K_3$ is well-known. The following theorem is surely known to many experts.

\begin{theorem}\label{B3class} \cite{BM}, \cite{TW} The algebra $K_3$ is semisimple over $\C$ except if the eigenvalues satisfy
one of the equations $\la_i=\la_j$, $\la_i^2-\la_i\la_j+\la_j^2=0$ or $\la_i^2+\la_j\la_k=0$, where $\{ i,j,k\} = \{ 1, 2, 3\}$.
In the semisimple case, we have three 1-dimensional representations $\{ \la_i\}$, three two-dimensional representations
$\{ \la_i\la_j\}$ and one three-dimensional representation $\{ \la_1\la_2\la_3\}$. Here, the expression between the
brackets $\{ \ \}$ is the determinant of $\sigma_1$ in the given representation. The central element $\Delta_3^2$ acts on these representations via the scalars $\la_i^6$, $-\la_i^3\la_j^3$ and $(\la_1\la_2\la_3)^2$.
\end{theorem}

$Proof.$ See e.g. \cite{TW} where the irreducible representations of $B_3$ of dimension $\leq 5$ are classified.

\subsection{Generically irreducible representations of $K_4$}\label{genericirrep} We call a representation of an algebra which algebraically depends on
parameters $x_1, x_2,\ ...\ x_d$ {\it generically irreducible},  if it is irreducible for a Zariski open subset
in the parameter space. The generically irreducible representations of $K_4$ have been classified, see \cite{MHecke}.
 To describe them, first observe that
the images of $\sigma_1$ and $\sigma_3$ form a commutative subalgebra of the image of $B_4$ in any representation
of $B_4$. This motivates the following definition.

\begin{definition}\label{weights} Let $W$ be a representation of $B_4$. We call any common eigenvector of $\sigma_1$
and $\sigma_3$ a {\it weight vector}. If $\sigma_1$ has eigenvalues $\la_1,\la_2,\ ...\ \la_d$, we say that a vector $v\in W$ has
weight $(i,j)$ if $\sigma_1$ and $\sigma_3$ act on it via multiplication by $\la_i$ and $\la_j$ respectively.  The dimension of the corresponding eigenspace
is called the multiplicity of the weight. 
\end{definition}

\begin{remark}\label{ijrem} It follows from Eq. \ref{braidflip} that $\Delta_4$ maps the weight space of $(i,j)$ onto the weight space of $(j,i)$.
\end{remark}
The following table lists the generically irreducible representations of $K_4$. As for $K_3$, each simple
representation is already characterized
by the determinant of a standard generator, except for the 9-dimensional
representations. Here, it also depends on the choice of a primitive
third root of unity $\theta$.  We will refer to the representation by $\{ {\rm det}(\sigma_1)\}$,
except for the nine-dimensional representation $\{ \la_1^3\la_2^3\la_3^3\}_\theta$ where we add the primitive
third root of unity.
  We list the scalar by which 
 $\Delta_4^2$ acts in the third  column below.
The weights listed in the fourth column all have multiplicity 1 except for $(1,1)$ in the 8-dimensional representation, 
which has multiplicity 2 as indicated. The last column lists the simple $B_3$ modules which appear in a composition series of simple modules if we restrict the action to $B_3$. This is taken from \cite{MWa}. We did not list the additional representations which are obtained via permuting the eigenvalues.
\[
\begin{tabular}{|c|c|c|c|c|} dim & det($\sigma_1$)&$\Delta_4^2$&weights& $B_3$ rep\\
\hline
& & & &\\
9 & $\la_1^3\la_2^3\la_3^3$&$ \theta\la_1^4\la_2^4\la_3^4
$& $(i,j),\ 1\leq i,j\leq 3$& see $R9$\\
8&$\la_1^4\la_2^2\la_3^2$&
$\la_1^6\la_2^3\la_3^3$&  (1,1): 2, $(i,j), i\neq j$: 1 & see $R8$\\
6&$\la_1^3\la_2^2\la_3$&
$\la_1^6\la_2^4\la_3^2$& $(1,i),(i,1)$ and $(2,2)$& see $R6$ \\
3 & $\la_1\la_2\la_3$ & $\la_1^4\la_2^4\la_3^4$& $(i,i)$&$\{\la_1\la_2\la_3\}$\\
3 & $\la_1^2\la_2$ & $\la_1^8\la_2^4$ & $(1,1),(1,2),(2,1)$&$\{\la_1\la_2\}+\{\la_1\}$\\
2 &  $\la_1\la_2$ & $\la_1^6\la_2^6$ & $(1,1),(2,2)$&$\{\la_1\la_2\}$\\
1 & $\la_1$ & $\la_1^{12}$ & $(1,1)$&$\{\la_1\}$\\
\end{tabular}
\]
\vskip .2cm
$R9:$\quad $\{\la_1\la_2\la_3\}+ \{\la_1\la_2\}+ \{\la_1\la_3\} + \{\la_2\la_3\}$,

$R8:$\quad  $\{\la_1\la_2\la_3\} + \{\la_1\la_2\} + \{\la_1\la_3\} + \{\la_1\}$,

$R6:$\quad  $\{\la_1\la_2\la_3\} +\{\la_1\la_2\} + \{\la_1\}$.
\vskip .2cm
\begin{remark} The representations in which a generator only has one or two eigenvalues appear in the representation
theory of the Iwahori-Hecke algebra of type $A_3$ and are well-known. The additional three- and six-dimensional representations
appear as summands in the so-called $BMW$ algebras, while the eight-dimensional one appears in a two-parameter algebra
connected with Lie types $E_N$, see \cite{Wexc}. Explicit matrix representations  for the 6, 8 and 9-dimensional 
representations can be found in \cite{MWa}, Table 4. 
\end{remark}

\subsection{Non-generic representations}\label{possiblereps}
The representations listed in the last subsection may have non-trivial subquotient representations when the eigenvalues satisfy certain equations. These equations were determined in \cite{MW3} and are listed in the table below, where $\{ i,j,k\}=\{ 1,2,3\}$ and $\theta$ is a primitive third root of unity.

\[
\begin{tabular}{|c|c|} representation & equations\\
\hline
& \\
$\{\la_i\la_j\}$ & $\la_i^2-\la_i\la_j+\la_j^2$\\
$\{\la_i^2\la_j\}$ & $\la_i^2+\la_j^2$ \\
$\{ \la_1\la_2\la_3\}$ & $\la_i^2+\la_j\la_k$ \\
$\{\la_i^3\la_j^2\la_k\}$&$\la_i+\la_k,\ \la_j+\la_k,\ \la_j^2+\la_i\la_k,\ \la_i^3-\la_j^2\la_k$\\
$\{\la_i^4\la_j^2\la_k^2\}$&  $\la_i^2\la_j-\la_k^3,\ \la_i^2\la_k-\la_j^3, \ \la_i^2-\theta^{\pm 1}\la_j\la_k$\\
$\{\la_1^3\la_2^3\la_3^3\}_\theta $ &$\la_i +\theta^{\pm 1}\la_j,\ \la_i^2-\theta\la_j\la_k$\\
\end{tabular}
\]
We now list the $B_4$ representations which can be obtained as subquotients of the ones listed in the previous subsections, provided they are not isomorphic to one of the regular representations.
As the weights $(i,j)$ and $(j,i)$ have the same multiplicity by Remark \ref{ijrem},
we only list one of them. Obviously, we also obtain representations
of the same dimensions if we permute the eigenvalues, which are not listed separately here.
The bar in the column for the determinant for the 3-dimensional
representation only indicates that the eigenvalue $\la_i$ behaves differently in this representation
than the other two eigenvalues. Thus, in the notation $\{\la_j|\la_i\la_k\}$, the distinguished eigenvalue is $\la_j$. This can be seen by looking at the weights.
The following table lists these new simple representations. The last column lists the simple factors in a composition series if we view the modules as a $B_3$ module.

\[
\begin{tabular}{|c|c|c|c|c|c|} equation& dim & notation&$\Delta_4^2$&weights & $B_3$\ {\rm rep}\\
\hline
& & & &\\
$\la_i^2=-\la_j^2$& 2 & $\{\la_i\la_j\}^*$&$\la_i^8\la_j^4$&$(i,j)$&$\{\la_i\la_j\}$\\
$\la_i=-\la_k$&3&$\{\la_j|\la_i\la_k\}$&$\la_i^6\la_j^4\la_k^2$&
 $ (j,j), (i,k)$ &$\{\la_1\la_2\la_3\}$   \\
$\la_j=-\la_k$&4&\{$\la_i^2\la_j\la_k\}$&$\la_i^6\la_j^4\la_k^2$&
 $(i,j), (i,k)$& $\{\la_1\la_2\la_3\},\ \{\la_i\}$\\
$\la_i^3=\la_j^2\la_k$& 5 & $\{\la_j^2|\la_i^2\la_k\}$ & $\la_i^6\la_j^4\la_k^2$& $(j,j), (i,j), (i,k)$& see $R5$\\
$ \la_i^2=\theta\la_j\la_k$&7 & $\{\la_i^3\la_j^2\la_k^2\}$ & $\la_i^6\la_j^3\la_k^3$&$(i,i),(i,j),(i,k), (j,k)$& see $R7$\\
\end{tabular}
\]
where
\vskip .2cm
$R5:$\quad  $\{\la_1\la_2\la_3\},\ \{\la_i\la_j\}$,

$R7:$\quad  $\{\la_1\la_2\la_3\},\ \{\la_i\la_j\},\ \{\la_i\la_k\}$.

\vskip .2cm
The decomposition into $B_3$ modules can be deduced from the corresponding results in Section \ref{genericirrep}. As our given representation appears as a subquotient in more than one
generic representation (which can be identified by its $\Delta_4^2$ value, see the corresponding column above), a $B_3$ module must appear in the restriction of each of these generic $B_4$ modules to $B_3$. It turns out that the intersection
has the correct dimension in all cases.
\subsection{No more exceptional representations} We have listed new simple representations of $B_4$ when the eigenvalues satisfy
one irreducible polynomial. These representations will in general not stay simple if they also satisfy a second equation.
Nevertheless, we have the following result.

\begin{theorem}\label{nomoredecomposition} (\cite{MW3}, Theorem 4.7) The tables for generic and non-generic representations list all possible
simple representations of $K_4$ for which $\la_i\neq \la_j$, $1\leq i\neq j\leq 3$.
\end{theorem}

\def\et{{\bf e}}
\def\ft{{\bf f}}
\subsection{Temperley-Lieb algebras}\label{TLalgebras} We review some well-known facts of the Temperley-Lieb algebras. See e.g. \cite{GW}, \cite{Wseq} for more details. The Temperley-Lieb algebra $TL_n$ is given via generators 1 and  $\et_i$, $1\leq i <n$ and relations $\et_i^2=\et_i$, $\et_i\et_j=\et_j\et_i$ for $|i-j|>1$, and $\et_i\et_{i\pm 1}\et_i=\frac{1}{[2]^2}\et_i$. Here $[2]=q+q^{-1}$ and $q\in \C-\{\pm i\}$. We list some well-known facts about TL algebras:

\begin{enumerate}
\item The algebras $TL_n$ are semisimple whenever $q^2$ is not a primitive $\ell$-th root of unity for $1< \ell\leq n$. 
\item If $TL_n$ is semisimple, its simple modules $W_\mu$ are labeled by the Young diagrams $\mu=[\mu_1,\mu_2]$ with $|\mu|=\mu_1+\mu_2=n$.
\item $W_\mu\cong W_{[\mu_1-1,\mu_2]}\oplus W_{[\mu_1,\mu_2-1]}$ as a $TL_{n-1}$ module.
\item If $\ft_n$ is the central idempotent for the summand $[n,0]$, we have $(\ft_n\et_n\ft_n)^2=\frac{[n+1]}{[2][n]}\ft_n\et_n\ft_n$.
\item We obtain a representation of the braid group $B_n$ by mapping $\sigma_i\mapsto \la_11 +(\la_2-\la_1)\et_i$, whenever $-\la_1/\la_2=q^2$.
\end{enumerate}

\section{Determine eigenvalues of $c_{VV}$}

Throughout this section we assume $\Ca$ to be a tensor category of type $G_2$ or $G_{2,k}$, with $k\geq 1$ and $k\neq  3,6$. We denote the eigenvalues via which  the braiding morphism $c_{V,V}$ acts on the subobjects
$V_{2\La_1}$, $V_{\La_2}$, $V_{\La_1}$ and $V_0=\1$ of $V^{\otimes 2}$ by $\la_i$, $1\leq i\leq 4$, in that order.

\subsection{The new part $V^{\otimes n}_{new}$} We call a simple submodule $V_\mu\subset V^{\otimes n}$  {\it new} if it has not appeared in any tensor power $V^{\otimes m}$ with $m<n$. It follows from our tensor product rules that in this case $|\mu|=n$, or equivalently, $\mu=a\La_1+b\La_2$ with $a+2b=n$. If we have a decomposition of $V^{\otimes n}$ into a direct sum of simple modules, we denote by $V^{\otimes n}_{new}$ the submodule which contains all the new modules. Then it follows by induction on $n$, using our restriction rule \ref{braidrestrict}, that $\sigma_i$ only acts via
eigenvalues $\la_1$ and $\la_2$ on $V^{\otimes n}_{new}$, for all $n\in \N$ for which  $V^{\otimes n}_{new}\neq 0$.
\begin{proposition}\label{TLnew} Let  $\Ca$  be a tensor category of type $G_2$ or $G_{2,k}$, and let $q^2=-\la_1/\la_2$. Then $\la_1\neq\la_2$, and we obtain a surjective map $TL_n\to \End(V^{\otimes n})_{new}$ provided that $q^2$ is not an $\ell$-th root of unity for $1<\ell\leq n$.
\end{proposition}

$Proof.$ Since $V^{\otimes 2}$ is multiplicity-free,
$c_{V,V}$ is diagonalizable.
We first assume $k\neq9$ and treat $k=9$ at the end
of the proof.
If $\la_1=\la_2$, $\sigma_i$ would act as $\la_11$ on $\Hom(V_{2\La_1+\La_2},V^{\otimes 4})$, for $i=1,2,3$.
In this case $\Delta_4$ would act as $\la_1^61$ on that module, i.e. it would have trace $3\la_1^6$.
This would contradict the fact that $\Delta_4$ permutes the images of the projections $p_{[2]}\otimes p_{[1^2]}$ and $p_{[1^2]}\otimes p_{[2]}$,
see Remark \ref{ijrem}, and hence could only have trace $\pm \la_1^6$. Hence $\la_1\neq \la_2$. 
In particular, we can define nontrivial idempotents $e_i$ as the eigenprojections for the eigenvalue $\la_2$ of the images of $\sigma_i$ in $\End(V^{\otimes n})$.

If the 2-dimensional $B_3$-module $\Hom(V_{[2,1]},V^{\otimes 3})$ were isomorphic to $\{\la_1\}\oplus \{\la_2\}$, our restriction rules \ref{braidrestrict} would imply $\Hom(V_{2\La_1+\La_2},V^{\otimes 4})\cong 2\{\la_1\}\oplus\{\la_2\}$. This would again contradict  Remark \ref{ijrem}. Hence 
 $\Hom(V_{[2,1]},V^{\otimes 3})\cong\{\la_1\la_2\}$, which is simple unless $-\la_1/\la_2=\theta$, a primitive 3rd root of unity. 
Also observe that $\End(V^{\otimes 3}_{new})$ only has two simple summands, labeled by $[3,0]$ and $[2,1]$.
Both $e_i$'s act as zero in the former, and as rank 1 idempotents in the latter one. Hence $e_1e_2e_1$ is a scalar multiple of $e_1$, which is equal to 
$\frac{1}{-\la_1/\la_2+2-\la_2/\la_1}=\frac{1}{[2]^2}$, by \cite{TW}; recall that $q^2=-\la_1/\la_2$.

We can now prove the surjectivity statement by induction on $n$ as follows: We just showed it for $n=2$ and $n=3$ in the previous paragraphs. For the induction step $n-1\to n$, we observe that $\Hom(V_\mu,V^{\otimes n})$
is a semisimple $TL_n$ module for any $\mu$ with $|\mu|=n$, according to our assumption for $q^2$. It decomposes as a $TL_{n-1}$ module into the direct sum $\Hom(V_{[\mu_1-1,\mu_2]},V^{\otimes n-1})\oplus \Hom(V_{[\mu_1,\mu_2-1]},V^{\otimes n-1})$.
By the induction assumption, the module 
$\Hom(V_{[\mu_1-1,\mu_2]},V^{\otimes n-1})$ is the irreducible
$TL_{n-1}$ module labeled by $[\mu_1-1,\mu_2]$, while
$\Hom(V_{[\mu_1,\mu_2-1]},V^{\otimes n-1})$ is the irreducible
$TL_{n-1}$ module labeled by $[\mu_1,\mu_2-1]$. But then it follows from the restriction rules for $TL_n$, see Section \ref{TLalgebras} (3), that $\Hom(V_\mu,V^{\otimes n})$ must be isomorphic to the simple $TL_n$ module labeled by $\mu$.

Now let $k=9$. The simple objects of categories of type $G_{2,9}$ are labeled by Young diagrams with $\leq 3$ boxes, see Section \ref{fusion:sec}. So it does not contain a simple object labeled by $2\La_1+\La_2$ corresponding to the Young diagram $[3,1]$, and we cannot use the argument above.

It follows from the definition of $\Delta_n^2$ that we have $\det(\pi(\Delta_n^2))=\det(\pi(\sigma_1)^{n(n-1)})$ for any finite-dimensional representation $\pi$ of $B_n$. By taking determinants, it follows from Eq. \ref{twistbraid} and \ref{scalarmatch} that 
\begin{equation}\label{determinants}
\Theta_\mu^{d(\mu,n)}=\det(\sigma_1)^{n(n-1)}\Theta_V^{nd(\mu,n)}
\end{equation}
where $d(\mu,n)=\dim \Hom(V_\mu,V^{\otimes n})$ and $\det(\sigma_1)$ is the determinant of $\sigma_1$ in the representation of $B_n$ on $\Hom(V_\mu,V^{\otimes n})$. Applying Eq. \ref{twistbraid} for $n=2$,
we obtain
\begin{equation}\label{2twist}
\Theta_V=\la_3^{-2},\quad \la_4^2=\Theta_V^{-2},\quad \Theta_{2\La_1}=\la_1^2\Theta_V^2,\quad \Theta_{\La_2}=\la_2^2\Theta_V^2.
\end{equation}
We can therefore express all $\Theta_\mu$ in terms of the eigenvalues for $|\mu|\leq 2$. Plugging these values into Eq. \ref{determinants} for $n=3$,
we obtain the equations
$$\la_2^6=1,\quad \la_2^4=\la_1^6\la_3^2,\quad \la_1^6\la_2^6\la_3^2=1.$$
If $\la_1^2=\la_2^2$, we obtain $\la_3^2=1$ from the first and last equations. We then also obtain $\la_1^2=1=\la_2^2$ from the second equation,
and $\la_4^2=1$ from Eq. \ref{twistbraid}, see above.
This would imply that $\Ca$ is symmetric.
But this is impossible: the fusion rules
give $FPdim(V)=(3+\sqrt{21})/2$, the positive root of $x^2=3x+3$, while symmetric fusion categories have integral
Frobenius-Perron dimensions
\cite{EGNO}, Corollary 9.9.11. Hence $\la_1\neq \la_2$ also for $k=9$. 

Suppose $N=\Hom(V_{[2,1]},V^{\otimes3})$ is not simple. Its submodule and quotient of dimension one
have braid eigenvalues $\la_1$ and $\la_2$.
Since the full twist acts by one scalar on $N$, we have $\la_1^6=\la_2^6=1$. The determinant equations above then give
$\la_3^2=1$ and $\la_2^4=1$, hence $\la_2^2=1$.
By Eq. \ref{2twist}, we also have
$\Theta_V=1$ and $\la_4^2=1$. Moreover, $\Delta_3^2$ acts as the identity on $N$,
so Eq. \ref{twistbraid} gives $\Theta_{[2,1]}=1$.
We must have $\la_1^2\neq1$, as otherwise all four
eigenvalues would have square one and $\Ca$
would be symmetric, which was ruled out above.

Let us first assume $\la_2=\la_3$.
The space $M=\Hom(V_{2\La_1},V^{\otimes3})$
has dimension three, and each of $\sigma_1$
and $\sigma_2$ is diagonalizable with eigenvalues
$(\la_1,\la_2,\la_2)$.
Their eigenspaces for $\la_2$ have dimension two
and hence have a nonzero intersection.
On a common eigenvector, $\Delta_3^2$ acts by
$\la_2^6=1$.
On the other hand, Eqs. \ref{twistbraid}
and \ref{2twist} show that it acts on all of $M$
via the scalar
$\Theta_{2\La_1}\Theta_V^{-3}=\la_1^2$.
Hence $\la_1^2=1$, a contradiction.

It remains to consider $\la_2=-\la_3$.
As $\la_2^2=\la_3^2=1$ and $\la_1^2\neq1$,
the eigenvalues $\la_1,\la_2,\la_3$ are distinct.
By Proposition \ref{8dimidentification},
$\Delta_4^2$ acts on
$W=\Hom(V_{[2,1]},V^{\otimes4})$
via the scalar
$$
\la_1^6\la_2^3\la_3^3=-1.
$$
But Eq. \ref{twistbraid} gives its scalar on $W$
as $\Theta_{[2,1]}\Theta_V^{-4}=1$.
This contradiction proves that $N$ is simple.
The same calculation with $TL_n$ above now proves
the claim for $n\leq3$.
For $n\geq4$, the new part is zero, since every
simple object has already appeared in a tensor
power of degree at most three.

\begin{lemma}\label{TLrootofunity}  Assume $q^2=-\la_1/\la_2$ is a primitive $\ell$-th root of unity.
Then $\Ca$ can not be of type $G_{2,k}$ if that type allows a simple object labeled by $[\ell,1]$.
\end{lemma}

$Proof$. By Proposition \ref{TLnew}, $TL_{\ell-1}$ maps surjectively onto $ \End(V^{\otimes \ell-1}_{new})$,
with the image $f_{\ell-1}$ of $\ft_{\ell-1}$ mapping onto $V_{[\ell-1,0]}\subset V^{\otimes \ell-1}_{new}$.
By our fusion rules, $f_{\ell-1}$ acts as a rank 1 idempotent in both $\Hom(V_{[\ell,0]},V^{\otimes \ell})$
and $\Hom(V_{[\ell-1,1]},V^{\otimes \ell})$. By Section \ref{TLalgebras},(4), $f_{\ell-1}e_{\ell-1}f_{\ell-1}$ is nilpotent, and hence acts as 0 in both modules. It therefore is in the kernel of the homomorphism from $TL_\ell\to \End(V^{\otimes \ell})_{new}$, and hence also in the kernel of the homomorphism from $TL_{\ell+1}\to \End(V^{\otimes \ell+1})_{new}$. But then also
$$e_\ell f_{\ell-1}e_{\ell-1}f_{\ell-1} e_\ell\ =\ \frac{1}{[2]^2} f_{\ell-1}e_\ell$$
is in the kernel. But by our tensor product rules, if $V_{[\ell,1]}$ is an object in $\Ca$, it is a summand in the tensor product $V_{[\ell-1,0]}\otimes V_{[1,1]}$, i.e. $f_{\ell-1}e_\ell$ would act nonzero in
$\Hom(V_{[\ell,1]},V^{\otimes \ell+1})$. This is a contradiction.

\begin{proposition}\label{Casimirprop} Let $\Ca$ be a tensor category of type $G_2$ or $G_{2,k}$
such that $TL_n\to \End(V^{\otimes n})_{new}$ is surjective for all $n$. Then the twisting numbers $\Theta^\Ca_\mu$ are already determined by $\Theta_V$ and the eigenvalues $\la_1$ and $\la_2$ of $c_{V,V}$. In particular, if $\la_1=q^2$, $\la_2=-1$ and $\Theta_V=q^{12}$, then we have $\Theta^\Ca_\mu=q^{C_\mu}$, the twisting number  of the corresponding simple object in $\U_q$ or $\bar{\U}_q$.
\end{proposition}

$Proof$. By Eq \ref{twistbraid}, it suffices to determine $\Theta_V$ and the scalar via which $\Delta_n^2$ acts on $V_\mu$ when it occurs for the first time, in $V^{\otimes |\mu|}$. 
By Proposition \ref{TLnew}, the scalar of $\Delta_n^2$ only depends on the representation theory of the Temperley-Lieb algebra $TL_n$. If the values of $\la_1$, $\la_2$ and $\Theta_V$ coincide with the ones
in $\U_q$, we therefore also get $\Theta^\Ca_\mu = \Theta^{\U_q}_\mu=q^{C_\mu}$, see Proposition \ref{prop:eigenvalue}.

\subsection{Properties of the $B_4$ module $\Hom(V_{\La_1+\La_2}, V^{\otimes 4})$}\label{Wproperties}
 Let $\Ca$ be of type $G_2$ or $G_{2,k}$, $k\neq 3,6$, and let $W=\Hom(V_{\La_1+\La_2}, V^{\otimes 4})$.
Then we have the following simple lemma, which is a consequence of the axioms for a ribbon tensor category.

\begin{lemma}\label{8dimprop} The representation of $B_4$ on $W$ satisfies the following properties:

(a) The full twist $\Delta_4^2$ acts via a scalar $\al$ on $W$.

(b) The full twist $\Delta_3^2$ acts via a scalar $\beta$ on each $B_3$ submodule $\cong (\Hom(V_{2\La_1}, V^{\otimes 3}))$ of $W$, see Eq \ref{braidrestrict},

(c) The module $W$ has the same weights as the generic 8-dimensional irreducible representation, see the table in Section \ref{genericirrep}.
\end{lemma}

\begin{remark}\label{exoticex} Properties (a)-(c) in Lemma \ref{8dimprop} are satisfied by any 8-dimensional representation which has a composition series with the factors of a composition series of the canonical 8-dimensional representation. Nontrivial composition series can occur when the eigenvalues satisfy certain relations. They were determined in \cite{MW3}, and are listed below for the reader's convenience.
\end{remark}
\[
\begin{tabular}{|c|c|c|} case&  relations&  simple factors\\
\hline&&\\
(1)&$(\la_2^3-\la_1^2\la_3)$&$\{\la_1^2|\la_2^2\la_3\}+\{\la_1^2\la_3\}$\\
(2)&$(\la_1^2-\theta\la_2\la_3)$&$\{\la_1^3\la_2^2\la_3^2\}+\{\la_1\}$\\
(3)&$(\la_2^3-\la_1^2\la_3, \la_3^3-\la_1^2\la_2,\la_2^2+\la_3^2)$&$\{\la_1^2\la_2\}+\{\la_1^2\la_3\}+\{\la_2\la_3\}^*$\\
(4)&$(\la_2^3-\la_1^2\la_3, \la_3^3-\la_1^2\la_2,\la_2+\la_3)$&$\{\la_1|\la_2\la_3\}+\{\la_1\la_2\}^*+\{\la_1\la_3\}^*+\{\la_1\}$\\
(5)&$(\la_2^3-\la_1^2\la_3, \la_1^2-\theta\la_2\la_3, \la_1^2+\la_3^2)$&$ \{\la_1^2|\la_2^2\la_3\}+\{\la_1\la_3\}^*+\{\la_1\}$\\
(6)&$(\la_2^3-\la_1^2\la_3, \la_1^2-\theta\la_2\la_3,\la_1+\la_3)$&$ \{\la_2^2\la_1\la_3\}+\{\la_1^2\la_3\}+\{\la_1\}$\\
\end{tabular}
\]
The cases obtained by interchanging the indices $2$ and $3$ are also allowed.
\begin{proposition}\label{8dimidentification} Assume that $\Ca$
is a tensor category of type $G_2$ or $G_{2,k}$, $k\geq 1$, $k\neq 3,6$ such that the eigenvalues $\la_1$, $\la_2$ and $\la_3$
are mutually distinct.  Then the $B_4$ module $W=\Hom(V_{\La_1+\La_2},V^{\otimes 4})$ 
satisfies the conditions of Lemma \ref{8dimprop} with $\al=\la_1^6\la_2^3\la_3^3$ and $\beta_{2\La_1}=\la_1^2\la_2^2\la_3^2$.
\end{proposition}

\subsection{Proof of Proposition \ref{8dimidentification}}\label{rulingout} We will use our classification of simple $K_4$ representations to prove Proposition \ref{8dimidentification}. This requires to consider
several cases.

\begin{lemma}\label{lessthanthree}
Let $\Ca$ be a tensor category of type $G_2$ or $G_{2,k}$. Assume that $W=\Hom(V_{\La_1+\La_2},V^{\otimes 4})$ satisfies the conditions of Lemma \ref{8dimprop} with $\la_i\neq\la_j$ for $i\neq j$. Then it  must either be the irreducible module $\{\la_1^4\la_2^2\la_3^2\}$, or it must allow a composition series with factors as listed in the previous table.
\end{lemma}

$Proof.$ Let us first assume that $W$ allows a composition series where none of its factors has dimension $>3$. Observe that the only simple $K_4$ representations of dimension $\leq 3$ which contain weights $(i,j)$
with $i\neq j$ are $\{\la_i^2\la_j\}$, $\{\la_i\la_j\}^*$ and $\{\la_k|\la_i\la_j\}$. Representations of dimension $\leq 3$
that contain weight $(1,1)$ but not $(2,2)$ or $(3,3)$ are $\{\la_1^2\la_i\}$, $i=2,3$, $\{\la_1|\la_2\la_3\}$ or $\{\la_1\}$. In view of these few cases, one can directly check that the only
possible cases would be

(a) Cases (3) and (4) in the table above,

(b) a composition series with factors $\{\la_1|\la_2\la_3\}+\{\la_1^2\la_2\} +\{ \la_1\la_3\}^*$, or a similar one with $\la_2$ and $\la_3$ interchanged.

Observe that in case (b) $\{\la_1|\la_2\la_3\}$ and $\{\la_1\la_3\}^*$ exist only if 
$\la_2=-\la_3$ and $\la_1^2=-\la_3^2$. This also implies that $\la_1^2=-\la_2^2$. But then the representation $\{\la_1^2\la_2\}$ is not simple, i.e. we are in case (4).

Hence, we can assume that a composition series of $W$ must contain at least one factor with dimension $\geq 4$. 
By our classification of representations of $K_4$, there is  only one simple 8-dimensional representation satisfying the 
conditions of Lemma \ref{8dimprop}. Similarly, we have shown in Section \ref{possiblereps} that the only 
 simple 7-dimensional representation of $B_4$ containing the weight (1,1) is
 $\{\la_1^3\la_2^2\la_3^2\}$. By the conditions for the weights of $W$, the remaining simple factor must be $\{\la_1\}$, which gives us case (2) in our table.
Any simple 6-dimensional representation contains weights of the form $(i,i)$ and $(j,j)$ for $i\neq j$. Hence, it can not be a factor in a composition
series of $W$. A 5-dimensional simple representation contains the weight $(j,j)$ for some $j$, $1\leq j\leq 3$.
Hence, it can appear in a composition series of $W$ only if $j=1$, i.e. it is isomorphic to $\{\la_1^2|\la_2^2\la_3\}$. If it appears as a factor in $W$, the remaining factor(s) must have weights $(1,1)$, $(1,3)$ and $(3,1)$. 
It follows from the tables in Section \ref{possiblereps}
that the only possible $K_4$ representations with these weights are $\{\la_1^2\la_3\}$ and $\{\la_1\}+\{\la_1\la_3\}^*$, i.e. we are in cases (1) or (5).

Assume now that $W$ contains a factor isomorphic to the 4-dimensional representation $\{\la_i^2\la_j\la_k\}$. Recall that such
representations only exist when $\la_j=-\la_k$. We need to consider different cases:

(1) If $\la_2=-\la_3$, the weights of our 4-dimensional factor are $(1,2),(2,1), (1,3)(3,1)$. Hence, the remaining factors need to have the weights $(1,1)$ (twice),
$(2,3)$ and $(3,2)$. Inspecting the possible simple factors for such a representation, there are two possibilities:

(a) We have the 1-dimensional representation $\{\la_1\}$ and the 3-dimensional representation $\{ \la_1|\la_2\la_3\}$. 
As $\Delta_4^2$ has to act by the same scalar on each factor, we obtain the equalities
$$\la_1^6\la_2^4\la_3^2=\la_1^6\la_2^2\la_3^4=\la_1^{12}=\la_1^4\la_2^6\la_3^2=\la_1^4\la_2^2\la_3^6.$$
As $\la_2=-\la_3$, we obtain $\la_1^6\la_2^6=\la_1^4\la_2^8$, from which we deduce $\la_1^2=\la_2^2$ and hence $\la_1=\pm \la_2$. 
But this would imply $\la_1\in \{\la_2,\la_3\}$.

(b) We assume that we have two factors isomorphic to $\{\la_1\}$ and one factor isomorphic to the exceptional 2-dimensional 
representation $\{\la_2\la_3\}^*$ with weights $(2,3)$ and $(3,2)$. The latter is possible only if $\la_2^2=-\la_3^2$, contradicting
our assumption $\la_2=-\la_3$.

(2) If $\la_1=-\la_3$, our 4-dimensional factor would have the weights $(1,2),(2,1), (2,3)(3,2)$, and the remaining factors
would have the weights $(1,1)$ (twice),
$(1,3)$ and $(3,1)$. We can rule out factors isomorphic to $\{\la_1\}$ (twice) and $\{\la_1\la_3\}^*$ as in the previous case 1(b).
The only other possibility would be factors isomorphic to $\{\la_1\}$ and $\{\la_1^2\la_3\}$, i.e. we are in case (6) of our table.

(3) If $\la_1=-\la_2$, we can use the proof of Case (2) by just interchanging $\la_2$ with $\la_3$.

\medskip
$Proof$ of Proposition \ref{8dimidentification}:  We have checked all possible 8-dimensional representations of $B_4$ which satisfy 
the conditions in Lemma \ref{8dimprop}. In all cases, $\Delta_4^2$ acts via the required scalar on $W$. Moreover, in all these cases, the restriction
to $B_3$ contains a 3-dimensional subrepresentation on which $\Delta_3^2$ acts via the required scalar $\la_1^2\la_2^2\la_3^2$. 
This follows from the restriction rules stated in Sections \ref{genericirrep} and \ref{possiblereps}. This finishes the proof of  Proposition \ref{8dimidentification}.

\begin{corollary}\label{3dimirrep} If $W$ satisfies the conditions of Lemma \ref{8dimprop}, it contains the irreducible $K_3$ module $\{\la_1\la_2\la_3\}$ in all cases except case (3) in the table before Proposition  \ref{8dimidentification}.
\end{corollary}

$Proof.$ This follows from the tables in Sections \ref{genericirrep}, \ref{possiblereps} and \ref{Wproperties}.

\subsection{Getting eigenvalues from the ribbon structure} This has already mostly appeared in \cite{MW2}, where we had assumed that the image of $B_3$ generates $\End(V^{\otimes 3})$. This told us via which scalars $\Delta_3^2$ has to act on the simple components of $V^{\otimes 3}$. Similar techniques also work for other exceptional Lie types, see \cite{MST}, Section 5.7, using results in \cite{TW}.

\begin{proposition}\label{eigenvaluesribbon}
Let $\Ca$ be a ribbon tensor category of type $G_2$ or $G_{2,k}$, $k\neq 3,6$. Assume that the eigenvalues of $\sigma_1$ in $\End(V^{\otimes 2})$
are given by $\la_1$, $\la_2$,  $\la_3$, $\la_4$ for $V_\mu\subset V^{\otimes 2}$ with $\mu=2\La_1, \La_2, \La_1$ and $0$
in that order.  
Then setting $\la_1=q^2$, we have $\la_2=-1$, $\la_3=- q^{-6}$ and $\la_4= q^{-12}$.
\end{proposition}

$Proof.$  Let $V_\mu\subset V^{\otimes n}$ and let $\Theta_\mu$ be the scalar of the twist map on $V_\mu$. Observe that Eqs \ref{determinants} 
and \ref{2twist} hold for arbitrary $k>0$, $k\neq 3,6$, from which we obtain
\begin{equation}\label{eqtw2}
\Theta_V=\Theta_{\La_1}=\la_3^{-2}\quad {\rm and}\quad \la_4^2=1/\Theta_V^2 = \la_3^4.
\end{equation}
We now assume that the eigenvalues $\la_1,\la_2,\la_3$ are mutually distinct. Then  Proposition \ref{8dimidentification} gives us the scalar via which $\Delta_3^2$ acts on $\Hom(V_{2\La_1},V^{\otimes 3})$. Applying Eq. \ref{scalarmatch} for $n=2$ and $n=3$ and Eq. \ref{2twist}, we obtain
$$\la_1^2\la_3^{-4}=\Theta_{2\La_1}=(\la_1\la_2\la_3)^2\la_3^{-6}.$$
It follows from this  that
\begin{equation}\label{eqtw4}
\la_2^2=1.
\end{equation}
Observe that $\Delta_3^2$ acts via the scalar $-(\la_1\la_2)^3$ on $\Hom(V_{\La_1+\La_2},V^{\otimes 3})$, by Theorem \ref{B3class} and  Proposition \ref{TLnew}. Using the same argument for $\Theta_{\La_1+\La_2}$ as before and  Proposition \ref{8dimidentification}, we obtain
\begin{equation}\label{eqtw3}
-(\la_1\la_2)^3\la_3^{-6}=\Theta_{\La_1+\La_2}=\la_1^6\la_2^3\la_3^{3-8},
\end{equation}
from which we derive
\begin{equation}\label{eqtw5}
\la_3=-\la_1^{-3}.
\end{equation}
Setting $\la_1=q^2$, we obtain $\la_3=- q^{-6}$ from Eq. \ref{eqtw5} and  $\la_4^2=q^{-24}$ from Eq. \ref{eqtw2}. 
The signs of $\la_2$ and $\la_4$ can now be determined as in the proof of \cite{MW2}, Proposition 3.2. The remaining cases where the first three eigenvalues of $c_{V,V}$ are not distinct will be treated in the next lemma.

\begin{lemma}\label{repeat} Let $\Ca$ be a category of type $G_2$ or $G_{2,k}$ such that $|\{ \la_1,\la_2,\la_3\}|\leq 2$, i.e. at most two of the eigenvalues are distinct.
Then $\Ca$ is symmetric. In particular, if $q^2$ is as in Proposition \ref{eigenvaluesribbon}, it can not be a primitive 3rd or 8th root of unity.
\end{lemma}

$Proof.$ (a) It follows from Proposition \ref{TLnew} that $\la_1\neq \la_2$ and that $\Theta_{\La_1+\La_2}=(-\la_1\la_2)^3\Theta_V^3$. In particular, $\la_3$ has to be equal to $\la_1$ or $\la_2$.

(b) By our assumption, the 3-dimensional representation $\Hom(V_{2\La_1},V^{\otimes 3})$ must be reducible.
Its simple composition factors therefore have
dimensions $1+1+1$ or $1+2$. As $\Delta_3^2$ has to act as a scalar on it, we conclude that $\la_1^6=\la_2^6$
or $(-\la_1\la_2)^3=\la_1^6$ or $(-\la_1\la_2)^3=\la_2^6$. In all of these cases, we have $\la_1^6=\la_2^6$ and $\Theta_{2\La_1}=\la_1^6\Theta_V^3=\la_2^6\Theta_V^3$.

(c) By our tensor product rules, the object $V_{\La_1+\La_2}$ appears in $V_{2\La_1}\otimes V_{\La_2}$. As $\la_1\neq \la_2$, it follows that the 8-dimensional representation $W=\Hom(V_{\La_1+\La_2},V^{\otimes 4})$
must contain weights $(1,2)$ and $(2,1)$. The only simple representations of $B_4$ containing such weights such that the braid generators only have eigenvalues $\la_1$ and $\la_2$ are the 3-dimensional representations $\{\la_1^2\la_2\}$ and $\{\la_2^2\la_1\}$. As $\dim\ W=8$, it must also contain either 1-dimensional or 2-dimensional simple representations. As $\Delta_4^2$ has to act via a scalar on $W$,
we deduce that $\la_1^8=\la_2^8$ by comparing all possible values of $\Delta_4^2$ for simple representations only depending on $\la_1$ and $\la_2$. We conclude from this and (b) that $\la_1^2=\la_2^2=\la_3^2$.

(d) As $\Theta_V=\la_3^2\Theta_V^2$, we obtain from (b) that $\Theta_V=1/\la_i^2$, $1\leq i\leq 3$.
On the other hand, we have $\Theta_{2\La_1}=\la_1^2\Theta_V^2=\la_1^2\la_1^{-4}=\la_1^{-2}$,
while we obtained in (b) that $\Theta_{2\La_1}=\la_1^6\Theta_V^3=1$. Hence $\la_i^2=1$ for $1\leq i\leq 3$. 
Finally, we also obtain $1=\Theta_0=\la_4^2/\la_3^4=\la_4^2$, which implies that our braiding is symmetric.

For the last statement, observe that $|\{q^2,-1,-q^{-6}\}|<3$, but  $\{q^4,1,q^{-12}\}\neq\{1\}$ for $q^2$ a primitive 3rd or 8th root of unity.

%\begin{corollary}\label{signchange} If $\al=\pm 1$ and $\la_1=-\la_2$, then $\{ \la_3, \la_4\}\subset  \{\la_1,\la_2\} =\{\pm 1\}$. \end{corollary}

%$Proof.$ We have $q^2=\la_1=-\la_2=1$. Hence also $\{\la_3,\la_4\}=\{ -\al q^{-6},\al^2 q^{-12}\}\subset \{\pm 1\}$.

\section{Getting dimensions from eigenvalues}

\subsection{Main result of section} We determine for a given tensor category of type $G_2$ or $G_{2,k}$ all possible values of $q^2$ and the corresponding dimension functions. The result is as follows.

\begin{theorem}\label{dimfromeigen} (a) If $\Ca$ is a ribbon category of type $G_2$, the eigenvalue $\la_1=q^2$ is either 1 or not a root of unity.

(b) If $\Ca$ is of type $G_{2,k}$, $k\geq 1$, $q^2$ must be a primitive $(k+12)$-th root of unity.
In particular, $\Ca$ is not symmetric.

(c) The categorical dimension of the simple object $V_\mu$ is equal to $d_\mu$ as in \ref{dimensionformula} for any $\mu\in P_{+}$ (case (a)) or $\mu\in P_{+,k}$ (case (b)).
\end{theorem}
The theorem will be proved for most values of the parameters in Proposition \ref{dimfromeigenprop}.
The remaining cases will be dealt with in the following sections. We already mention the following immediate corollary.

\begin{corollary}\label{TLss:cor}
The representations $TL_n\to \End(V^{\otimes n})_{new}$ in Proposition \ref{TLnew} are always surjective.
\end{corollary}

$Proof.$ This follows from Proposition \ref{TLnew} for $\Ca$ of type $G_2$, as $TL_n$ will always be semisimple for $q$ as in part (a) of the theorem. If $\Ca$ is of type $G_{2,k}$, a $TL_n$ module $W_\la$ appears in  $V^{\otimes n}_{new}$ only if $n = |\la|\leq k+6<k+12=\ell$. Hence again the image of $TL_n$ is semisimple for $q^2$ a primitive $\ell$-th root of unity, see Section \ref{TLalgebras}.

%$Proof$. This was proved in Proposition \ref{dimfromeigenprop} under the restrictions stated there. The cases with two of the eigenvalues of $c_{V,V}$ being the same was treated in Proposition \ref{repeat} and Section \ref{repeated2}. This already ruled out $q^2$ being a third or 6th root of unity. We can not use the argument in the proof of Proposition \ref{dimfromeigenprop}  if $q^2$ is a fourth, fifth, or ninth root of unity, as the denominator in the formula for $d_\mu$ would be zero in that case. These cases were ruled out in Sections \ref{repeated2} and \ref{ruleout4}.
\subsection{The generic case}\label{generic:sec} The following proposition determines the dimensions for most tensor categories of type $G_2$ or $G_{2,k}$.

\begin{proposition}\label{dimfromeigenprop} The statement in Theorem \ref{dimfromeigen} holds if the eigenvalues of $c_{V,V}$, as stated in Proposition \ref{eigenvaluesribbon}  are mutually distinct, and $q^2$ is not an $\ell$-th root of unity, $\ell\in \{ 3,4,5,6,9\}$.
\end{proposition}

$Proof.$ This is a generalization of the proof of \cite{MW2}, Theorem 3.3.
It was shown in \cite{TW}, Section 3.2 that the eigenvalues of $c_{V,V}$ for a self-dual object $V$ and $\Theta_V$
determine the dimensions of the simple objects in $V^{\otimes 2}$, as long as there are at most 5 direct summands,
which are mutually non-isomorphic such that the braiding morphism $c_{V,V}$ acts via distinct eigenvalues on these summands. By Proposition \ref{eigenvaluesribbon}, the eigenvalues of the braiding morphism $c_{V,V}$ coincide with the ones of the corresponding
morphism in $\U_q$. This proves the claim in particular for the fundamental objects $V_{\La_1}$ and $V_{\La_2}$.
As these objects generate the representation ring for $\U_q$ or $\bar{\U}_q$, it follows that the dimension of a simple object $V_\mu$ is given by $d_\mu$ as defined in \ref{dimensionformula}. This proves statement (c) of Theorem \ref{dimfromeigen}.

In order to prove the restrictions for the values of $q$, we first recall an elementary direct proof of the fact that the dimension of any simple object is determined by the dimensions of the fundamental objects $V_{\La_i}$, $i=1,2$. We proceed  by induction on $|\mu|=\mu_1+\mu_2$. By the previous discussion, we have already shown the claim for $|\mu|\leq 2$.
Assume now $|\mu|=n+1\geq 3$ with $\mu_2\geq 1$. Then also $\nu=\mu-\La_2$ is a dominant integral weight. It follows from \ref{tensoradjoint} that 
\begin{equation}\label{tensorad}
V_\nu\otimes V_{\La_2}\ \cong\ V_\mu\oplus \bigoplus_\ga V_\ga,
\end{equation}
with $|\ga|\leq n$. By induction assumption, $\dim V_\ga = d_\ga$, as defined in \ref{dimensionformula}, for $|\ga|\leq n$.
Hence 
$$\dim V_\mu=\dim V_\nu \dim V_{\La_2}-\sum_\ga \dim V_\ga\ =\ d_\nu d_{\La_2}-\sum_\ga d_\ga\ =\ d_\mu,$$
by the tensor product rules for quantum groups; here the summation over $\ga$ is as in \ref{tensoradjoint}. If $\mu=(n+1,0,-n-1)=(n+1)\La_1$,
we can prove the claim by considering $V_{n\La_1}\otimes V$ and using the already proved dimension formula for $V_\ga$ with $\ga_2\geq 1$ or $|\ga|\leq n-1$.

Let now $\Ca$ be of type $G_{2,k}$ with $3|k$, and let $n=k/3$. We can show that
 $\dim V_\mu=d_\mu$ for $|\mu|\leq n$ using the proof in the last paragraph. Let $|\mu|=n+1$ with $\mu_2\geq 1$, and $\nu=\mu-\La_2$. By the tensor product rules for categories of type $G_{2,k}$, the right hand side of \ref{tensorad} does not contain the object $V_\mu$ anymore. Hence $d_\mu=d_\nu d_{\La_2}-\sum d_\ga=0$. It follows from Corollary \ref{dimestimate} that all factors in the numerator of $d_\mu$ are nonzero except possibly $[3(\mu_1+\mu_2+3)]=[k+12]$. This proves the claim if $3|k$.

Let now $\Ca$ be of type $G_{2,k}$ with $3\nmid k$, and $k$ odd. It follows from the definition of $P_{+,k}$ that all weights $\mu$ with $|\mu|\leq (k+7)/2$ are in $P_{+,k}$ except $[(k+7)/2,0]$. Using \ref{tensorrule} with $\la=[(k+5)/2,0]$,
one shows that $d_{[(k+7)/2,0]}=0$ and $[k+12]=0$ by a similar argument as in the previous paragraph.

If $\Ca$  is of type $G_{2,k}$ with $3\nmid k$, and $k$ even, all diagrams $\la$ with $|\la|\leq k/2+4$
are in $P_{+,k}$ except $[k/2+4,0]$ and $[k/2+3,1]$.
Using \ref{tensorad} with $\nu=[k/2+2,0]$, one shows as before that $d_{[k/2+3,1]}=0$. Again, using Corollary \ref{dimestimate}, one deduces that $[k+12]=0$ and $q^2$ is a primitive $(k+12)$-th root of unity.

\subsection{$q^2$ a primitive 4-th, 5-th or 9-th root of unity}\label{459sec} We have to be more careful if $q^2$ is a primitive $\ell$-th root of unity with $\ell\in \{ 3,4,5,6,9\}$, for which the denominator of the dimension formula becomes zero. We need the following elementary lemma.

\begin{lemma}\label{459}
Let $q^2$ be a primitive $\ell$-th root of unity, and let $d_\mu$ be as in Eq \ref{dimensionformula}. Then we have

(a) If $\ell=4$, $d_\mu\neq 0$ for all $\mu=[\mu_1,\mu_2]$ with $|\mu|\leq 3$ except for $\mu=[2,1]$.

(b) If $\ell=9$, $d_\mu\neq 0$ for all $\mu=[\mu_1,\mu_2]$ with $|\mu|\leq 6$ except for $\mu=[4,2]$.

(c)  If $\ell=5$, $d_\mu\neq 0$ for all $\mu=[\mu_1,\mu_2]$ with $|\mu|\leq 12$ except for $\mu=[8,4]$.
\end{lemma}

$Proof.$ Let $m_1=\mu_1+2$ and $m_2=\mu_2+1$. Observe that $m_1>m_2$ for any dominant weight $\mu$. Then the numerator in the formula for $d_\mu$ is equal to 
$$[m_1-m_2][2m_1+m_2][m_1+2m_2][3m_1][3m_2][3m_1+3m_2].$$
In order to prove (c), we must check that only one factor in the numerator becomes zero for a primitive 5-th root of unity for all weights $\mu$ with $|\mu|\leq 12$, except for $\mu=[8,4]$.
In the latter case, it is easy to check that each factor in the numerator will become zero for $q^2$ a primitive 5-th root of unity, hence $d_{[8,4]}=0$. In the former case, it is probably easiest to list for each factor $[f]$ all Young diagrams $\mu$ with $|\mu|\leq 12$ for which $5|f$. One then observes that each such Young diagram only appears once in these six lists.

The cases with $\ell=4$ and $\ell=9$ can be checked similarly and are easier.

\begin{corollary}\label{no459}
We do not obtain a tensor category of type $G_2$ or $G_{2,k}$ if $q^2$ is a primitive fifth or ninth root of unity. The same is true for $q^2$ a primitive fourth root of unity if the dimensions for the fundamental objects are as in Lemma \ref{459}.
\end{corollary}

$Proof.$ One checks easily that the eigenvalues of $c_{V,V}$ are distinct if $q^2$ is a primitive fifth or ninth root of unity. Using this and our assumption for $\ell=4$, we can assume $\dim V_\mu=d_\mu$ for all simple objects $V_\mu$ in a possible category $\Ca$ of type $G_{2,k}$. Lemma \ref{459} would then imply that we have objects labeled by
$[3]$, $[5,1]$ and $[9,3]$ for $\ell=4,9,5$ respectively, but, by rigidity, no objects labeled by $[2,1]$, $[4,2]$ and $[8,4]$. But it follows from the definitions  that $\mu=[\mu_1,\mu_2]\in P_{+,k}$ always implies $[\mu_1-1,\mu_2+1]\in P_{+,k}$, a contradiction to the last sentence.

\begin{remark}\label{459remark} (a) As the four eigenvalues of $c_{V,V}$ are not distinct for $q^2$ a primitive 4-th root of unity, it might be conceivable that there exists a tensor category for our given braiding with different dimensions. This will be ruled out in Section \ref{ruleout4}.

(b) One can obtain interesting tensor categories from braid representations for roots of unity as in Lemma \ref{459}. E.g. for $\ell=4$, one can construct a non-symmetric tensor category with the fusion ring of the symmetric group $S_4$. We plan to study such categories in \cite{RW}.
\end{remark}

\subsection{Repeated eigenvalues}\label{repeated2}  We have already shown in Lemma \ref{repeat} that there are no non-symmetric tensor categories of type $G_2$ or $G_{2,k}$ with $\la_i=\la_j$, $1\leq i\neq j\leq 3$,
which also ruled out that $q^2$ is a primitive third root of unity. 
Hence it remains to prove Theorem \ref{dimfromeigen} when $\la_4\in \{\la_1,\la_2,\la_3\}$, which also covers the cases with $q^2$ a primitive $\ell$-th root of unity for $\ell=4,6$.

\begin{lemma}\label{repsplit}
Under the assumptions of Proposition
\ref{eigenvaluesribbon}, suppose that
$\la_4\in\{\la_1,\la_2,\la_3\}$ and that
$\la_1,\la_2,\la_3$ are distinct.
Then $q^4=-1$, and
$$
\Hom(V,V^{\otimes3})
\cong \{\la_2\}\oplus\{\la_1\la_2\la_3\},
$$
where the two summands are simple of dimensions one
and three, respectively.
\end{lemma}

$Proof.$ Let $t=q^2$ and
$W=\Hom(V,V^{\otimes3})$. By Proposition \ref{eigenvaluesribbon}, the eigenvalues
of $\sigma_1$ on $W$ are
$t,-1,-t^{-3},t^{-6}$.
The equality of the last eigenvalue with one of the
first three gives
$t^7=1$, $t^6=-1$, or $t^3=-1$.
Since the first three eigenvalues are distinct, $t$
has order $4$, $6$, $7$, or $12$. By Eq. \ref{twistbraid}, $\Delta_3^2$ acts on $W$
as $\Theta_V^{-2}=t^{-12}$.
Since $\sigma_1$ is diagonalizable, the corollary to
\cite{TW}, Proposition 2.2, shows that every simple
factor of $W$ has dimension at most three and has
distinct eigenvalues. At least one such factor contains the eigenvalue $t$.
By \cite{TW}, Proposition 2.5, its full-twist scalar
is one of $t^6$, $t^3$, $t^{-6}$, $t^{-4}$,
corresponding to the eigenvalue lists
$(t)$, $(t,-1)$, $(t,-t^{-3})$, and
$(t,-1,-t^{-3})$, respectively.
Equality with $t^{-12}$ would require $t^{18}=1$, $t^{15}=1$, $t^6=1$, or $t^8=1$, respectively.
None is possible when $t$ has order $7$ or $12$. If $t$ has order $6$, then $\Theta_V=t^6=1$ and
the eigenvalues of $c_{V,V}$ are $t,-1,1,1$.
The ribbon trace identity, applied to both the braiding and the  reverse
braiding, gives
$$
0=Tr(c_{V,V}-c_{V,V}^{-1})
=(t-t^{-1})\dim V_{2\La_1}; 
$$
see \cite{EGNO}, Exercise 8.10.15.
This contradicts Lemma \ref{nonzero trace}, since
$t\neq t^{-1}$. Thus $t$ has order $4$. Theorem \ref{B3class} shows
that $K_3(t,-1,-t^{-3})$ is semisimple. The only simple modules on which $\Delta_3^2$ acts
as $t^{-12}=1$ are the one dimensional module
$\{\la_2\}$ and the three dimensional module
$\{\la_1\la_2\la_3\}$. Their eigenvalue multiplicities give the stated decomposition.
\\

\begin{lemma}\label{exotic eigenvalues}
Assume that $\Ca$ is a tensor category of type $G_{2,k}$ such that  $\la_4$ coincides with one of the eigenvalues $\la_i$, $i<4$. Then $\Ca$ is a symmetric tensor category, or $q^4=-1$, with $q$ as in Proposition \ref{eigenvaluesribbon}. In particular, $q^2$ cannot be a primitive sixth, seventh, or 12th root of unity.
\end{lemma}
$Proof.$ If the first three eigenvalues are not
distinct, Lemma \ref{repeat} shows that $\Ca$ is
symmetric and $q^4=1$. Otherwise, Lemma \ref{repsplit}
gives $q^4=-1$. In either case, $q^2$ cannot have
order $6$, $7$, or $12$.\\

\subsection{The case with eigenvalues $(\la_1,\la_2,\la_3,\la_4)=(\sqrt{-1},-1,-\sqrt{-1},-1)$}\label{ruleout4} We now deal with the case mentioned in Remark \ref{459remark}. For simplicity of notation, we fix $i=\sqrt{-1}$. 

\begin{lemma}\label{noexoticdim} Let $\Ca$ be a ribbon tensor category of type $G_2$ or $G_{2,k}$ with an object $V$ whose braiding and tensor product rules for $V^{\otimes 2}$ and $V^{\otimes 3}$ are as in Proposition \ref{eigenvaluesribbon}, with $q^2=i$. Then the dimensions for $\mu$ with $|\mu|\leq 3$ are given by Eq \ref{dimensionformula}.
\end{lemma}

$Proof.$ It is straightforward to check, using Theorem \ref{B3class} that $K_3(i,-1,-i)$ is semisimple. By Lemma \ref{repsplit}, the $B_3$ module $W=\Hom(V,V^{\otimes 3})$ splits into the direct sum of a simple 3-dimensional module $\{ i,-1,-i\}$ and a one-dimensional module $\{-1\}$. It will be convenient to use the labeling $\ga_1=i$, $\gamma_2=-i$, and $\ga_3=\ga_4=-1$. We use bases of eigenvectors for both subrepresentations, starting with the 3-dimensional one. It follows from Theorem \ref{AB2theorem} (see \cite{MW3}, Section 3.6, Example (1) for more details) that the diagonal entries
 $d_i(\ga_3)$ of the eigenprojection of $S_2$ in the 3-dimensional summand are equal to
$$d_1(\ga_3)\ =\ d_2(\ga_3)\ = \frac{1}{2},\quad d_3(\ga_3)=0.$$
Let us denote the projection onto the trivial subobject $\1\subset V^{\otimes 2}$ by $e$.
Observe that $e$ is a sub-idempotent of the eigenprojection of $\ga_3=-1$.
It follows that $e_1$ has nonzero diagonal entries only for the 3rd and 4th basis vectors,
say $\al$ and $\beta$,
while $e_2$ only has nonzero diagonal entries for the first, second and fourth basis vectors.
As $e_2=\Delta_3 e_1\Delta_3^{-1}$, it follows that the last diagonal entry of $e_2$ is also $\beta$, and its first and second entries are equal to $\al/2$. We deduce from this
\begin{equation}\label{dimbeta}
\frac{1}{\dim^2 V}e_1\ =\ e_1e_2e_1\ =\ \beta^2e_1.
\end{equation}
On the other hand, as the projection onto $V\subset V^{\otimes 2}$ is equal to the eigenprojection for $\ga_2=-i$, we obtain
\begin{equation}\label{dimalpha}
\dim V\ =\ \frac{\al}{2}\dim^2 V\quad \Leftrightarrow\quad \frac{1}{\dim V}\ =\ \frac{\al}{2}.
\end{equation}
Since $e_1\neq 0$, Eq.~\ref{dimbeta} gives
$(\dim V)^{-2}=\beta^2$, and hence $\dim V=\pm 1/\beta$.
Moreover, $e_1$ is a rank-one idempotent whose only nonzero diagonal
entries are $\al$ and $\beta$. Therefore
$\operatorname{Tr}(e_1)=\al+\beta=1$.
If $\dim V=1/\beta$, it follows from $\al+\beta=1$ and the two previous equations that
$\al/2=\beta=1-\al$. This implies $\al=2/3$ and $\beta=1/3$.
We conclude
$$\dim V\ =\ 3\ =\ \dim V_{2\La_1},\quad \dim V_{\La_2}=2,$$
where the last equation follows from the decomposition of $V^{\otimes 2}$. 

If $\dim V=-1/\beta$, the same calculations give
$\dim V=\dim V_{2\La_1}=1$ and $\dim V_{\La_2}=-2$.
Since $q^2=i$, Eq. \ref{eqtw2} gives $\Theta_V=-1$.
The ribbon trace identity gives
$$
Tr(c_{V,V}^{-1})=\Theta_V^{-1}\dim V=-1;
$$
see \cite{EGNO}, Exercise 8.10.15.
On the other hand, the eigenvalues of $c_{V,V}^{-1}$
on $V_{2\La_1},V_{\La_2},V,\1$ are $-i,-1,i,-1$.
Thus
$$
Tr(c_{V,V}^{-1})
=-i\dim V_{2\La_1}-\dim V_{\La_2}+i\dim V-1
=1,
$$
a contradiction.
\medskip

{\it Proof of Theorem \ref{dimfromeigen}} We have shown that $q^2$ can not be an $\ell$-th root of unity for any tensor category of type $G_2$ or $G_{2,k}$, $k\geq 1$, if $\ell= 3,8$ (Lemma \ref{repeat}), $\ell=4,5,9$ (Corollary \ref{no459} and Lemma \ref{noexoticdim}) or $\ell=6, 7, 12$ (Lemma \ref{exotic eigenvalues}). 
If $q^2=\pm 1$, the tensor category $\Ca$ would be symmetric. This would entail integer dimensions in case of fusion tensor categories. This is not possible, see \cite{Ru}, Theorem 4.1.
This rules out all the cases which have not already been covered in the proof of Proposition \ref{dimfromeigenprop}.

\section{Surjectivity and uniqueness of braid representations}

It was shown in \cite{MW2} for tensor categories of type $G_2$ that surjectivity of the map $\C B_3\to \End(V^{\otimes 3})$ implies surjectivity of $\C B_n\to \End(V^{\otimes n})$ for all $n\in \N$. We will prove this result here also for any tensor category of type $G_{2,k}$, without any initial assumption for $n=3$. This will again be proved using quantum Jucys-Murphy type elements, extending methods already used in \cite{Wexc} and \cite{MW2}.

%\subsection{Representations with $A$ having at most 2 eigenvalues} 
\subsection{Representations of $AB_2$}\label{AB2section}
We now consider representations of $AB_2$ given by matrices $A$ and $T$ on a vector space $W$ which satisfy   the following conditions:
\vskip .2cm
\renewcommand{\labelenumi}{\alph{enumi})}
\begin{enumerate}
\item They satisfy the braid relation $ATAT=TATA$,
\item The matrix $A$ satisfies the equation $(A-\lambda_1I)(A-\lambda_2I)(A-\lambda_3I)=0$, where 
the eigenprojection $P$ for the eigenvalue $\lambda_3$ has rank 1.
\item We have $TATA=\delta I$.
\item We assume that $T$ is a diagonal matrix with eigenvalues $x_i$,
$1\leq i\leq n=\dim W$ with $x_i\neq x_j$ if $i\neq j$.
\end{enumerate}
It follows from these conditions that
\begin{equation}\label{APrelation}
A+\lambda_1\lambda_2A^{-1}=(\la_1+\la_2)I + \la_3^{-1}(\la_3-\la_1)(\la_3-\la_2)P.
\end{equation}
If \(p_{ij}\) denotes the \((i,j)\)-entry of the matrix \(P\), we deduce from the previous equation and Eq. \ref{Ainverse} that
\begin{equation}\label{linequa}
\frac{\delta + \la_1\la_2x_ix_j}{\delta} a_{ij} = (\la_1+\la_2)\delta_{ij} + \la_3^{-1}(\la_3-\la_1)(\la_3-\la_2)p_{ij}.
\end{equation}

The following theorem is a somewhat more precise reformulation of \cite{MW3}, Theorem 3.3 which goes back to \cite{Wexc}, Proposition 5.7.
\vskip .3cm
\begin{theorem}\label{AB2theorem} We assume a representation $W$ of $AB_2$ satisfying conditions (a)-(d) with $(\la_1-\la_3)(\la_2-\la_3)\neq 0$. 
If $\delta\not\in\{\la_k^2x_r^2,\ k=1,2,\ 1\leq r\leq n\}$ and  
$\delta+\la_1\la_2x_rx_s\neq 0$ for $1\leq r\neq s\leq n$, then

(a) the diagonal entries $d_r(\la_3)$ of $P$ are nonzero for all $1\leq r\leq n$,

(b) they are uniquely determined by
\begin{align}
d_r(\la_3)\ =\ &\frac{\la_3}{(\la_3-\la_1)(\la_3-\la_2)}\ \prod_{s\neq r}\frac{\delta+\la_1\la_2x_rx_s}{\delta(x_r-x_s)}\ \times \cr
&\times\ [(-1)^{n-1}\la_3\frac{\la_1\la_2x_r+\delta x_r^{-1}}{\delta}\ (\prod_{t=1}^n x_t)
- (\la_1+\la_2)(\frac{-\delta}{\la_1\la_2})^{(n-1-\ve)/2}x_r^\ve],
\end{align}
where $\ve=0$ for $n$ odd and $\ve=1$ for $n$ even. In particular, $W$ is an irreducible $AB_2$ module.
(c) The matrix $A$ is uniquely determined by the
diagonal entries of $P$, up to conjugation by a diagonal
matrix.

(d) The uniqueness assertions in (b) and (c) also
hold, with the other assumptions unchanged, if the
condition $\delta\neq\la_k^2x_r^2$, $k=1,2$, is replaced by
$p_{rr}\neq0$ for every $r$. 
\end{theorem}

$Proof.$ Parts (a) and (b) of this theorem were proved in \cite{Wexc}, Lemmas 4.6, 5.5 and 5.6, when $\la_1=q$ and $\la_2=-q^{-1}$. That proof also works in our setting.

For part (c), write $P=vw^t$. Since
$p_{rr}=v_rw_r\neq0$ for every $r$, all coordinates
of $v$ and $w$ are nonzero. Conjugating by
$D=\operatorname{diag}(v_1,\ldots,v_n)$ gives
$(D^{-1}PD)_{rs}=p_{ss}$. Thus the diagonal entries
determine $P$ up to diagonal conjugation. After this conjugation, Eq. \ref{linequa} determines
every off-diagonal entry of $A$. Since the vector
$(1,\ldots,1)^t$ is now a $\la_3$-eigenvector of $A$,
we have
$$
a_{rr}=\la_3-\sum_{s\neq r}a_{rs}.
$$ Thus the diagonal entries are determined as well.
This proves (c).

For part (d), write $P=vw^t$. Then all coordinates of $v$
and $w$ are nonzero. If
$\delta+\la_1\la_2x_r^2=0$, Eq. \ref{linequa}
determines $p_{rr}$ directly.
For the remaining indices, Eq. \ref{linequa} and
$Av=\la_3v$ give a linear system for the remaining
diagonal entries, with coefficient matrix
$$
\left(\frac{\delta}
{\delta+\la_1\la_2x_rx_s}\right)_{r,s},
$$
see \cite{Wexc}, Lemma 4.6.  This is an invertible Cauchy matrix, since the $x_r$
are distinct and its denominators are nonzero. 
The argument for (c) is then unchanged.

\subsection{Recent $AB_2$ representations}\label{pathexample} We assume that $\Ca$ is non-symmetric and,
in the fusion case, that $k\geq1$ and $k\neq3,6$. It will be instrumental to calculate the matrices $A=S(\la,\nu)$, see Lemma \ref{affinebraidlemma}, in the following setting:
Let $\la$ and $\nu$ be labels at level $n-1$ and $n+1$ in our $G_2$ or $G_{2,k}$ path diagram such that $\nu=\la+\om$ with $\om\in \{\ep_1-\ep_3=\La_1, \ep_2-\ep_3\}$ and with $|\nu|=n$ and $|\la|=n-1$; the notion $recent$ comes from the observation that $V_\nu$ has appeared in the $n$-th tensor power of $V$ for the first time. Then we obtain a representation of $AB_2$ on the vector space $W(\la,\nu)$ spanned by a basis labeled by the set $\P(\la,\nu)$ of paths of length 2 from $\la$ to $\nu$, see Sections \ref{pathrep} and \ref{affine:sec}. The rank of the eigenprojection $P$ of $A$ for the eigenvalue $\la_3$ corresponding to $V_{\La_1}\subset V^{\otimes 2}$ is equal to one, by Remark \ref{newpath}. To describe these representations, let $\om_1$, $\om_2$ be nonzero weights of $V$ such that $\om=\om_1+\om_2$. It follows from our description of weights that $\om_1$ and $\om_2$ span a parallelogram with diagonal $\om$ which splits it into two equilateral triangles. There are different cases which we first describe for type $G_2$.

(a) {\it Generic case} Assume that both $\la+\om_i\in P_+$ for $i=1,2$. 
Then  the paths
$r\in \P(\la,\la+\om)$ are of the form $r:\la\to\la+\phi_r\to\nu$, with
$$\phi_r\in I(\la,\la+\om)\ =\ \{ 0, \om, \om_1, \om_2\},$$
Observe that the multiplicity of nonzero weights in $V_{\La_2}$ is 1. It is not hard to check that the multiplicity of nonzero weights in $V_{2\La_1}$ is at most 2. Hence it follows that $A$ has three eigenvalues, $q^2$ with multiplicity 2, and $-1$ and $-q^{-6}$ with multiplicity 1 in this case, by Remark \ref{newpath}.

(b) It follows from the tensor product rules, see Section \ref{tensorG2}
that $\P(\la,\la+\om)$ has exactly three paths if and only if $\la=(a,0)$, with $a>0$ and $\om=\ep_1-\ep_3$. We leave it to the reader to check that $A$ acts with three different eigenvalues in these cases.

(c) The set $\P(\la,\la+\om)$ has exactly two paths if and only if either $\la=(a,a)$ with $a>0$ and $\om=\ep_1-\ep_3$ or $\la=(a+1,a)$ with $a\geq0$ and
$\om=\ep_2-\ep_3$. In both cases $A$ acts on the corresponding path space with eigenvalues $q^2$ and $-q^{-6}$, see the discussion in Section \ref{tensorG2}. 

(d) {\it Type $G_{2,k}$, $3|k$} It is straightforward to check that the same rules as above also apply for tensor categories of type $G_{2,k}$, if $3|k$, as long as both $\la$ and $\la+\om$ are in $P_{+,k}$. 

(e) {\it Type $G_{2,k}$, $3\nmid k$ } The same rules apply again if both $\la$ and $\la+\ep_1-\ep_3$ are in $P_{+,k}$,
except when $\om=\La_1$ and
$(\la+\rho,\La_1)=\ell-3$. Writing $\la=[a,b]$,
we have three paths with
eigenvalues $q^2,-1,-q^{-6}$, if $a>b>0$;
two paths with eigenvalues $-1,-q^{-6}$ if $a>0$
and $b=0$; and one path with eigenvalue $-q^{-6}$
if $a=b>0$. The eigenvalues follow from
Remark \ref{newpath}(a) and the fusion rules
in Section \ref{fusion:sec}. 

If  $(\la+\rho,\La_1)=\ell-2$ and $\om=\ep_2-\ep_3$,  we have 2 paths in $\P(\la,\la+\ep_2-\ep_3)$ with intermediate diagrams $\la$ and $\la-\ep_1+\ep_2$. Indeed, it is probably easiest to see geometrically  that $V_\nu\not\subset V_\la\otimes V_{\La_2}$ for $\nu=\la+\ep_2-\ep_3$, as the reflection in the boundary hyperplane cancels the weight $\la+\La_2+\rho$ with $\la+\ep_2-\ep_3 +\rho =\nu +\rho$; see e.g. the picture in  Section \ref{fusion:sec} for $\la=3\La_1$ and $\ell=13$. This also implies that $A$ only acts with eigenvalues $q^2$ and $-q^{-6}$.
If  $(\la+\rho,\La_1)=\ell-1$, $\la+\om\not\in P_{+,k}$ for our choices of $\om$.

\subsection{Checking conditions of Theorem \ref{AB2theorem}}
We assume the notations of the last subsection.

\begin{lemma}\label{divisibility} Assume $\la\in P_{+,k}$ such that also $\la+\psi\in P_{+,k}$ for all $\psi\in I(\la,\nu)$, and let $\ell=k+12$. Then we have

(a) $|(\la+\rho,\phi)|+2<\ell$ and $(\la+\rho,\phi)+1\neq 0$, $\phi\in I(\la,\nu)$

(b) $\ell\nmid (\la+\rho, \phi_1-\phi_2),\ \phi_1\neq\phi_2\in I(\la,\nu)$,

(c)  $\ell\nmid (\la+\rho, \phi_1-\phi_2)+1, \phi_1,\phi_2\in I(\la,\nu)$ if $3|\ell$
\end{lemma} 

$Proof.$ 
As $\phi\in I(\la,\nu)$ is a short root and $\La_1=\ep_1-\ep_3$ and $\La_2$ are the highest short and long roots, we have $\phi\leq \La_1<\La_2$ in Bruhat order. Moreover, we have $\ep_2-\ep_3=\La_1-(\ep_1-\ep_2)$. Hence $I(\la,\la+\om)$ contains $\La_1$ in both cases for $\om$.
As also $\la+\psi\in P_{+,k}$ for all $\psi\in I(\la,\nu)$, we have 
$$|(\la+\rho,\phi)|+2\leq (\la+\rho,\La_1)+2= (\la+\La_1+\rho,\La_1) <\ell,$$
as $\la+\La_1\in P_{+,k}$.

To prove (b), observe that 
$$\phi_1-\phi_2\in \{  \pm \om,\pm \om_1,\pm \om_2, \pm(\om_1-\om_2)\}.$$
Observe that these all are short roots except for $\pm(\om_1-\om_2)$. By definition of $P_{+,k}$, the claim follows if $\phi_1-\phi_2$ is a short root or if $3|\ell$. As $\om_1-\om_2$ is a long root, we have
%$$3|(\la+\rho,\om_1-\om_2)|\leq |(\la+\rho,\om_1)|+|(\la+\rho,\om_2)|<2\ell.$$
$$3\mid(\la+\rho,\om_1-\om_2),\qquad
|(\la+\rho,\om_1-\om_2)|
\leq |(\la+\rho,\om_1)|+|(\la+\rho,\om_2)|<2\ell.$$
But then $\ell\nmid(\la+\rho,\om_1-\om_2) $ if $3\nmid\ell$, which finishes the proof of (b).

Similarly, claim (c) follows from (a)  if $\phi_1-\phi_2$ is a short root or equal to $0$. If $3|\ell$, $(\la+\rho, \phi_1-\phi_2)+1$ would not be divisible by 3 if  $\phi_1-\phi_2$ is a long root.

\begin{lemma}\label{twiststructure} Assume $\Ca$ is a tensor category of type $G_2$ or $G_{2,k}$. Let $q^2$ be a primitive $\ell$-th root of unity, $\ell=k+12$, and $\la\in P_{+,k}$, where $k=\infty$ stands for $\Ca$ of type $G_2$ with $q$ not a root of unity. Then 

(a) $\Theta_{\la+\phi_1}\neq \Theta_{\la+\phi_2}$ for $\phi_1\neq\phi_2\in I(\la,\nu)$ such that also $\la+\phi_i\in P_{+,k}$ for $i=1,2$.

(b) Let $x_r=q^{c_{\la+\phi_r}-nc_V}$, $\phi_r\in I(\la,\nu)$. Then $\delta-q^2x_rx_s\neq 0$, except possibly for one $r=s\in\P(\la,\nu)$.

(c)  $\delta\neq \ga^2 x_r^2$ for $\ga\in\{ q^2,-1\}$ except possibly for one path $r$  if $\ga=-1$, $3\nmid \ell$, $\om=\ep_2-\ep_3$.
\end{lemma}

$Proof.$  It was shown in Proposition \ref{Casimirprop}  that $\Theta_\mu=q^{c_\mu}$ (we switch from $C_\mu$ to $c_\mu$, as there is no risk of confusing this with the braiding). 
It follows from the definition of $c_\la$ that
\begin{equation}\label{Casimirdiff}
c_{\la+\phi_1}-c_{\la+\phi_2}\ =\ 
\begin{cases} 2(\la+\rho, \phi_1-\phi_2) & {\rm if}\ \phi_i\neq 0,\ i=1,2,\cr
2(\la+\rho, \phi_1)+2& {\rm if}\ \phi_1\neq 0=\phi_2.
\end{cases}
\end{equation}
Claim (a) follows if $c_{\la+\phi_1}\not\equiv c_{\la+\phi_2}\ {\rm mod}\ 2\ell$.
This is a consequence of Lemma \ref{divisibility}, (a) if $\phi_2=0$, and Lemma \ref{divisibility}, (b) if $\phi_2\neq0$. 

%$\phi_1-\phi_2$ is a short root. If it is a long root and $3|\ell$, we also have  $$0<|(\la+\rho, \phi_1-\phi_2)|\leq (\la+\rho, \La_2)<\ell, $$  which implies the claim. Finally, if $3nmid|\ell$

It follows from Lemma \ref{affinebraidlemma} that $\delta=q^{c_\nu+c_\la-2nc_V}$. Let $\mu_i=\la+\phi_i$, $i=1,2$ be the intermediate weights for the paths $r$ and $s$. Recall that $\la_1\la_2=-q^2$. Then we have to show that $q^{c_\nu+c_\la-2nc_V}\neq q^2q^{c_{\mu_1}-nc_V}q^{c_{\mu_2}-nc_V}$, or, equivalently
$$2\ell\ \nmid \ c_\nu+c_\la-c_{\mu_1}-c_{\mu_2}-2=
2(\la+\rho,\om-\phi_1-\phi_2)+(\om,\om)-(\phi_1,\phi_1)-(\phi_2,\phi_2)-2.$$
The condition above is equivalent to
\begin{equation}\label{Casdif}
2\ell\ \nmid\ 
\begin{cases} 2(\la+\rho,\om) & {\rm if}\ \phi_1=\phi_2=0,\cr
2(\la+\rho,\om-\phi_1)-2 & {\rm if}\ \phi_1\neq 0=\phi_2,\cr
2(\la+\rho,\om-\phi_1-\phi_2)-4& {\rm if}\ \phi_2\neq 0.
\end{cases}
\end{equation}
This is obviously true in the first case.
One checks that $\om-\phi_1$ is a short root or 0 in the second case with $\phi_2=0$, and that $\om-\phi_1-\phi_2\in \{ 0, \pm\om_1,\pm\om_2, -\om,\}$ in the third case, if $\phi_1\neq \phi_2\neq 0$. So the claim again follows from 
Lemma \ref{divisibility}, (a) in these cases. We will do the case $r=s$ after part (c). 

For part (c), the claim for $\ga=-1$ and for $\ga=q^2$ is equivalent to 
$$2\ell \ \nmid \ c_\nu+c_\la-2c_\mu\hskip 2em {\rm and}\hskip 2em
2\ell \ \nmid \ c_\nu+c_\la-2c_\mu -4.
$$
This is translated to the first and third cases in Eq \ref{Casdif}, with adding/subtracting $ 2$ to the right hand side for $\ga=-1$/$\ga=q^2$ and $\phi_1=\phi_2=\phi$. Again, the case with $\phi=0$ follows from Lemma \ref{divisibility}, (a).

If $\om=\ep_1-\ep_3$, we have $\pm(\om_1-\om_2)=\pm (-\ep_1+2\ep_2-\ep_3)$  and if $\om=\ep_2-\ep_3$, we have $\pm(\om_1-\om_2)=\pm (2\ep_1-\ep_2-\ep_3)$. These are the only cases where $\om-2\phi$ is a long root. They are majorized by $\La_2$ in Bruhat order. It follows in this case  for $3|\ell$ that
$$|(\la+\rho, \om_1-\om_2)|\leq (\la+\rho, \La_2)-3<\ell-6,
$$
where the last inequality is derived from the fact that also $\la+\ep_1-\ep_3\in P_{+,k}$.
If $3\nmid\ell$ and $\om=\ep_1-\ep_3$, we use Lemma \ref{divisibility}, (a) to show that
$$0<(\la+\rho,-\ep_1+2\ep_2-\ep_3)= (\la+\rho,\ep_1-\ep_3)-2(\la+\rho,\ep_1-\ep_2)<\ell-3.$$
This and Lemma \ref{divisibility}, (a) proves the claim in this case.

If $3\nmid\ell$ and $\om=\ep_2-\ep_3$, the only possible long roots $\om-2\phi$ are $\pm(2\ep_1-\ep_2-\ep_3)$. If $\ga=q^2$, we have
$$
|\pm 2(\la+\rho,2\ep_1-\ep_2-\ep_3)-6|
=6(\la_1+2\mp1).
$$
Hence $2\ell$ would have to divide $2(\la_1+2\pm 1)$, as $3\nmid \ell$. But $0<\la_1+1<\la_1+3<\ell$,
by the definition of $P_{+,k}$. If $\ga=-1$, then $2\ell$ can divide at most one of $\pm 2(\la+\rho,2\ep_1-\ep_2-\ep_3)-2$, as their sum is $-4$, which is not divisible
by $2\ell$.

The missing case $r=s$ in part (b) is shown in the same way: we have the same set of roots $\om-2\phi$, where $2\ell$ can at most divide one of the two quantities $\pm 2(\la+\rho, \ga)-4$ for $\pm\ga$ is a long root of the form $\om-2\phi.$

\begin{lemma}\label{recent2:lem}Let $W$ be a recent $AB_2$ representation as in
Section \ref{pathexample}. Then the diagonal entries
of $P$ are nonzero and uniquely determined.
\end{lemma}
$Proof.$ If there is only one path, then $P=1$,
so the claim is immediate. Suppose that $p_{rr}=0$.
Since $P$ has rank one, its $r$-th row or column is
zero. If the column is zero, Eq. \ref{linequa} and
Lemma \ref{twiststructure}(b) show that the $r$-th coordinate line is preserved by $A$. The eigenvalue
$\ga$ on this line is not $\la_3$, since $P$ vanishes
on it. Equation \ref{Ainverse} thus gives
$\delta=\ga^2x_r^2$, with $\ga\in\{q^2,-1\}$.
If the row is zero, apply the same argument to
$A^t$ and $P^t$, using $T^t=T$. Hence  $\delta\neq \ga^2x_r^2$ for $\ga=\la_1,\la_2$ implies $p_{rr}\neq 0$. 

It follows from the last paragraph and Lemma \ref{twiststructure}(c) that $p_{rr}=0$ only if $3\nmid\ell$,
$\om=\ep_2-\ep_3$ and $\ga=-1$. It follows from our discussion in Section  \ref{pathexample} that this is not possible if we only have three or two paths (see also Lemma \ref{2by2matrixlem}).  For four paths, the quotient by the corresponding invariant line has eigenvalues $q^2, q^2, -q^{-6}$. By Lemma \ref{2by2matrixlem} this is a direct sum of representations of dimension at most two. Since $-q^{-6}$ has multiplicity one, there must be a one-dimensional subquotient with eigenvalue $q^2$. This would imply  $\delta=q^4x_s^2$ for some path $s$ contrary to Lemma \ref{twiststructure}(c).
Thus all diagonal entries of $P$ are nonzero. They are uniquely determined by $\la_1,\la_2,\la_3,\delta$ and the $x_r$, by
Theorem \ref{AB2theorem}(b),(d), together with
Lemma \ref{twiststructure}(a),(b).

Since the diagonal entries are nonzero, their
uniqueness follows from Theorem \ref{AB2theorem},(d)
together with Lemma \ref{twiststructure}.

\begin{theorem}\label{surjective:thm} Let $\Ca$ be a non-symmetric ribbon tensor category of type $G_2$ or $G_{2,k}$ with $k\geq 1$, and let $V$ be a simple object corresponding to the 7-dimensional representation of the Lie algebra $\g(G_2)$. Then the map $\C B_n\to \End(V^{\otimes n})$ is surjective for all $n\in \N$.
\end{theorem}

$Proof.$ This is proved by induction on $n$, see also the proof of  \cite{MW}, Theorem 2.7 (there is a misprint, indicating the end of the proof when it continues).  The claim is obviously true for $n=1$ and $n=2$, as we have ruled out all cases where we have fewer than four distinct eigenvalues for $c_{V,V}$.

For the induction step $n\to n+1$, first observe that $\Hom(V_\nu, V^{\otimes n+1})$ is a simple $B_{n+1}$ module for $|\nu|=n+1$, by Corollary \ref{TLss:cor}. Moreover, it follows from Lemma 
\ref{path:lem} and induction assumption that $\Hom(V_\nu, V^{\otimes n+1})$ is simple for $|\nu|\leq n-1$. We will prove the simplicity of  $\Hom(V_\nu, V^{\otimes n+1})$ for the remaining cases $|\nu|=n$ in a similar way, as outlined in
Remark \ref{newpath}.

Let us first assume $\nu_2\geq 1$ and $\nu_1> \nu_2+1$. Then $\Hom(V_\nu, V^{\otimes n+1})$
decomposes as a $B_n$ module into the direct sum of the simple modules allowed by the fusion rules,
with the possible labels $\nu, [\nu_1+1,\nu_2-1], [\nu_1-1,\nu_2+1], [\nu_1,\nu_2-1], [\nu_1-1,\nu_2]$. By the induction hypothesis, the map
$\C B_{n-1}\to\End(V^{\otimes n-1})$ is surjective.
We obtain two recent $AB_2$ representations, whose bases are labeled by paths of length two from $\la =[\nu_1-1,\nu_2]$ to $\nu$, and from  $\la=[\nu_1,\nu_2-1]$ to $\nu$, see Section \ref{pathexample}. By Lemma \ref{recent2:lem}, the braid matrix $A=S_n$ has an eigenprojection $P$ of rank 1 with nonzero diagonal entries in each of these $AB_2$ representations. This implies that the  $B_n$ submodules $\Hom(V_\mu, V^{\otimes n})$ corresponding to the intermediate diagrams of the remaining paths for a given $AB_2$ representation are in the same $B_{n+1}$ submodule of $\Hom(V_\nu, V^{\otimes n+1})$, see Remark \ref{newpath}. Together, the two sets of intermediate diagrams
label all the $B_n$ summands that occur. Both sets contain
$[\nu_1-1,\nu_2]$ and $[\nu_1,\nu_2-1]$, since the loops at these weights and the edge between them are not removed
at the affine boundary. Thus all the summands are connected,
and $\Hom(V_\nu,V^{\otimes n+1})$ is simple. 
%The statement about the diagonal entries now follows from Theorem \ref{AB2theorem}.

The proof in the cases where $\nu_2=0$ goes the same way. We just have fewer subdiagrams and paths. In the remaining cases $\nu=[a+1,a+1]$, $\la=[a+1,a]$ and $\nu=[a+1,a]$, $\la=[a,a]$, we only have two paths from $\la$ to $\nu$, with intermediate diagrams $[a+1,a]$ and $[a+2,a]$ in the first case. As already stated in the proof of Lemma \ref{recent2:lem} (see Section \ref{tensorG2} and Remark \ref{newpath}), $-1$ is not an eigenvalue of $A$ in this case. Hence $\delta\neq\ga^2x_r^2$ by Lemma \ref{twiststructure},(c), and $P$ and $A$ have nonzero diagonal entries, by Lemma \ref{2by2matrixlem}. Hence $\Hom(V_{(a+1,a+1)}, V^{\otimes 2a+3})$ and $\Hom(V_{(a+1,a)}, V^{\otimes 2a+2})$ are also simple $B_{2a+3}$ and $B_{2a+2}$-modules, respectively. 
% Hence a path idempotent corresponding to a fixed path ending at $\la$ is a minimal idempotent in the braid-group image, so Example \ref{pathexample} and Lemma \ref{affinebraidlemma}(b) apply.

\section{Classification of tensor categories of type $G_2$ or $G_{2,k}$}\label{classification1}  
\subsection{Uniqueness of braid representations} We can now proceed as in \cite{MW2} to show that the surjective braid representations in Theorem \ref{surjective:thm} for any tensor category of type $G_2$ are uniquely determined by the given fusion rules and the eigenvalue $\la_1=q^2$, up to isomorphism.  

\begin{theorem}\label{unique:thm} Let $\Ca$ be a non-symmetric ribbon tensor category of type $G_2$. Then the representations $B_n\to \End(V^{\otimes n})$ are uniquely determined by the eigenvalue $\la_1=q^2$ of $c_{V,V}$, acting on the summand $V_{2\La_1}\subset V^{\otimes 2}$.
\end{theorem}

$Proof.$ This theorem was proved in \cite{MW2}, Theorem 4.1 for tensor categories of type $G_2$, under the assumption that the map $\C B_3\to \End(V^{\otimes 3})$ is surjective. It was shown in Theorem \ref{surjective:thm} that this is always the case.

\medskip
For reasons which will be discussed in Remark \ref{general case:rem} we will prove the analogous result for ribbon tensor categories of type $G_{2,k}$ only for $n\leq 4$.

\begin{lemma}\label{noah} Let $\Ca$ be a ribbon tensor category of type $G_{2,k}$, $k\geq 1$. Then the representations $B_n\to\End(V^{\otimes n})$, $1\leq n\leq 4$, are uniquely determined by the value of $q^2$.
\end{lemma}

$Proof.$ We can assume that $q^2$ is a primitive $\ell$-th root of unity, $\ell=k+12$, by Theorem \ref{dimfromeigen}.
We have already shown that the eigenvalue $\la_1=q^2$ determines the eigenvalues of $c_{V,V}$, which proves the claim for $n=2$. By Theorem \ref{surjective:thm}, all nonzero $B_3$ modules $\Hom(V_\nu,V^{\otimes 3})$ are simple, with the action of $\Delta_3^2$ on $\Hom(V,V^{\otimes 3})$ equal to $\Theta_V^{-2}$, by Eq. \ref{twistbraid}. By \cite{TW}, Theorem 2.9 all simple representations of $B_3$ of dimension $\leq 5$ are uniquely determined by the eigenvalues of a braid generator and, for dimension $\geq 4$, also of $\Delta_3^2$. This proves the claim for $n=3$. 
For $n=4$, we can similarly use the classification of simple representations of $B_4$, provided the braid generator has at most three distinct eigenvalues, see Theorem \ref{nomoredecomposition}. This proves the uniqueness of the $B_4$ module $\Hom(V_\nu, V^{\otimes 4})$ for $|\nu|\geq 3$. 

It remains to prove the uniqueness of the $B_4$ module $\Hom(V_\nu, V^{\otimes 4})$ for $|\nu|\leq 2$. This will be done as in \cite{MW2}, see the outline in Section 4.2. The idea is to show that we can scale our path basis for the category $\Ca$ such that we get exactly the same matrices for the braid generators as in the corresponding path representations of the given tensor category $\U_q$ of type $G_2$ or $G_{2,k}$.

Let us do this for the 12-dimensional representation $\Hom(V_{[2,0]}, V^{\otimes 4})$, the most complicated case. 
By induction assumption, we can assume that the matrices for $\sigma_1$ and $\sigma_2$ coincide in both path representations. Let us consider the $3\times 3$ block $S_3([1,1],[2,0])$, acting on the paths with intermediate diagrams $\mu\in \{ [1,0], [2,0],[2,1]\}$; see Sections \ref{pathrep} and \ref{AB2section} for notation. It follows that $x_r\in \{ q^{12}, q^{28}, q^{42}\}$ and $\delta=q^{52}$. As usual, we set $\la_1=q^2$, $\la_2=-1$ and $\la_3=-q^{-6}$. It follows from our tensor product rules, see Section \ref{tensorG2} and Remark \ref{newpath} that each eigenprojection $P(\la_i)$ has rank 1. Hence we can use Theorem \ref{AB2theorem} to calculate $A$ via each of these eigenprojections, provided its assumptions are satisfied.
If $q^2$ is a primitive $\ell$-th root of unity with $\ell> 12$, one checks directly that 
$\delta\neq -\la_i\la_jx_rx_s$ for $r\neq s\in\P([1,1],[2,0])$.
Moreover, $\delta=\ga^2x_r^2$ for $\ga=-1$ only if $\ell=14,16$, for $\ga=q^2$ only if $\ell=18$ and for $\ga=-q^{-6}$ only if $\ell=20$. Hence we can find for each $\ell>12$ at least two eigenvalues $\la_j,\la_k$ such that $\delta\neq \ga^2x_r^2$ for $\ga\in \{\la_j,\la_k\}$ and $r\in \P([1,1],[2,0])$. It follows from Theorem \ref{AB2theorem} that the eigenprojection $P(\la_i)$, $\la_i\neq \la_j,\la_k$ has uniquely determined nonzero diagonal entries which determine the matrix block 
 $S_3([1,1],[2,0])$ up to conjugation by a diagonal matrix.
The analogous statement for the $3\times 3$ block $S_3([1,0],[2,0])$ can be shown similarly. 

We are now in the situation to prove uniqueness of our braid representation using the results in \cite{MW2}, Section 4.2. The cited lemmas use nothing more than checking the braid relation $\sigma_2\sigma_3\sigma_2=\sigma_3\sigma_2\sigma_3$. 

(1) We know the diagonal entries of the $3\times 3$ blocks of $S_3$ coincide. Using Lemma 4.6 in \cite{MW2}, we can make sure that both $3\times 3$ blocks as well as the $1\times 1$ blocks of $S_3$ coincide with the ones in representations for $\U$. The proof works for $3\times 3$ blocks as well as it was done for  $4\times 4$ blocks in \cite{MW2}.

(2) Using Lemmas 4.7 and 4.8 in \cite{MW2}, one can then show that also the $5\times 5$ block of $S_3$ can be made to coincide with the one in $\U$. Even though the proofs in \cite{MW2} were stated for $7\times 7$ blocks, they work as well in our setting.

The proofs of the uniqueness of the braid representations on $\Hom(V_\nu,V^{\otimes 4})$ for $\nu=[0,0]$, $\nu=[1,0]$ and $\nu=[1,1]$ go in a similar way and are easier.

\begin{remark}\label{general case:rem} It is quite straightforward to deduce our main result, Theorem \ref{classification-theorem} for categories of type $G_2$ or $G_{2,k}$ whenever Theorem \ref{unique:thm} holds, see Section \ref{main:sec}.
We have checked that one can prove it also for categories of type $G_{2,k}$. However, the proof would be even more tedious than the one of Lemma \ref{noah}.
Instead, we will give an alternative proof of our main result Theorem \ref{classification-theorem} in Section \ref{altunique}, using a different description of categories of type $G_2$, going back to Kuperberg \cite{Ku2}. It will only require  Lemma \ref{noah}.
\end{remark}
\subsection{Definition of functor}\label{functordef} We now show how our main theorem can  be derived from Theorems \ref{surjective:thm} and for whatever cases one can prove Theorem \ref{unique:thm}, essentially following the procedure in \cite{KW} (or \cite{TW}). As one has to make some minor adjustments, we give a fair amount of details about it. 

\def\Vt{{\tilde V}}
\def\At{{\tilde{\mathcal A}}}
\def\Xt{{\tilde X}}
\def\Yt{{\tilde Y}}
\def\it{{\tilde i}}
\def\dt{{\tilde d}}
\def\pt{{\tilde p}}
\def\vt{{\tilde v}}
\def\ut{{\tilde u}}
\def\tt{{\bf t}}
\def\at{{\bf a}}
\def\bt{{\bf b}}
\def\st{{\bf s}}
\def\sqt{{\bf sq}}
\def\ft{{\bf f}}
\def\gt{{\bf g}}
\def\pb{{\bf p}}

We have shown that for $\Ca$ of type $G_2$, and $V$ an object corresponding to the simple 7-dimensional representation of $\g(G_2)$, there exists a $q^2$ such that the representations of $\C B_n\to \End(V^{\otimes n})$ are isomorphic to those of $\C B_n\to \End(\Vt^{\otimes n})$, where $\Vt$ is the 7-dimensional simple object in $\U_q$. This means we have isomorphisms
$$\Phi: \A_n=\End_\Ca(V^{\otimes n})\to  \At_n=\End_{\U_q}(\Vt^{\otimes n}),$$
given by mapping the image of the $i$th braid generator in $\A_n$ to its image in $\At_n$, $1\leq i<n$.
We define the functor $F:\Ca\to\U_q$ as follows:

Given an object $X$ in $\Ca$, we define $F(X)=\Xt$, where $\Xt$ is an object in $\U_q$ whose image in the fusion ring is isomorphic to the one of $X$. If $f:X\to Y$ is a morphism in $\Ca$, we define $F(f)$ as follows: Pick $n$ large enough such that we can find monomorphisms $i_X:X\to V^{\otimes n}$ and $i_Y:Y\to V^{\otimes n}$ and epimorphisms $d_X:V^{\otimes n}\to X$ and $d_Y:V^{\otimes n}\to Y$. We normalize them
such that $p_X=i_Xd_X$ and $p_Y=i_Yd_Y$ are idempotents in $\A_n$.

Let $\Xt=F(X)$ and $\Yt=F(Y)$. We similarly define morphisms $\it_X,\it_Y,\dt_X,\dt_Y$ such that
$\it_X\dt_X=\pt_X=\Phi(p_X)\in \End(\Vt^{\otimes n})$ and $\it_Y\dt_Y=\pt_Y=\Phi(p_Y)$. We fix the inclusion and projection maps for
each object and use the compatible extensions to
larger tensor powers described below. The morphism $F(f)$ is then defined by
\begin{equation}\label{Fmorph}
F(f)\ =\ \dt_Y\Phi(i_Yfd_X)\it_X.
\end{equation}

\begin{lemma}\label{Fwell-defined} The construction above defines a functor $F$.
Different choices of $n$ and of the splittings give
naturally isomorphic functors.
\end{lemma}
$Proof.$ For $f:X\to Y$ and $g:Y\to Z$,
choose a common tensor power containing $X,Y,Z$.
Equation \ref{Fmorph} gives $F(g)F(f)=F(gf)$ and
$F(1_X)=1_{F(X)}$. Thus $F$ is a functor. Let $F'$ use another set of choices
in the same tensor power. For each object $X$, put
$$
\eta_X=\dt'_X\Phi(i'_Xd_X)\it_X.
$$
Using $d_Xi_X=d'_Xi'_X=1_X$,
$\dt_X\it_X=1_{F(X)}$ and
$\dt'_X\it'_X=1_{F'(X)}$, we see that
$\eta_X$ is invertible, with inverse
$$
\eta_X^{-1}=\dt_X\Phi(i_Xd'_X)\it'_X.
$$
For every $f:X\to Y$, Eq. \ref{Fmorph} gives
$$
\eta_YF(f)
=\dt'_Y\Phi(i'_Yfd_X)\it_X
=F'(f)\eta_X.
$$
Hence the two functors are naturally isomorphic. 

To prove the general statement, consider morphisms $i_k: \1\to V^{\otimes k}$ and $d_k:V^{\otimes k}\to \1$, for $k=2,3$. All these morphisms are unique up to  scalar multiples. We assume as before that $d_ki_k=1$, for $k=2,3$. We leave it to the reader to check that replacing $i_X$ by $i_X\otimes i_k$, and $d_X$ by $d_X\otimes d_k$ for the definition of $F(f)$ (with similar adjustments for $i_Y$ and $d_Y$) leads to a naturally isomorphic functor, for $k=2,3$. For these extensions, we also replace $\it_X$
by $\it_X\otimes\it_k$ and $\dt_X$ by
$\dt_X\otimes\dt_k$, and similarly for $Y$, where
$\dt_k\it_k=1$ and $\it_k\dt_k=\Phi(i_kd_k)$.
Since $\Phi$ respects tensor products, these compatible
extensions leave $F(f)$ unchanged.
Other choices in the enlarged tensor power are
compared by the maps $\eta_X$ above.
Any two tensor powers can be enlarged to a common
one, since every integer at least two is a sum
of twos and threes. This, together with the result proved in the last paragraph, implies the claim.

\subsection{Symmetric categories}\label{symmetric:sec} We have already ruled out symmetric tensor categories of type $G_{2,k}$ for $k\geq 1$ in Theorem \ref{dimfromeigen}. For type $G_2$, it is known to experts that there is only one symmetric tensor category of that type, up to equivalence. We indicate a proof for the readers' convenience below, which was outlined to us by Pavel Etingof. Possible errors are due to the authors.

\begin{proposition}\label{sym:prop} There is only one symmetric semisimple tensor category of type $G_2$, up to equivalence.
\end{proposition}

$Proof.$  Given a tensor category $\Ca$ of type $G_2$
with object $V$ corresponding to the 7-dimensional representation of the group $G$ of type $G_2$, it follows that
$$\dim \End_\Ca(V^{\otimes n})\leq 7^{2n}.$$
Hence we can apply Deligne's theorem \cite{D1}, \cite{Ostrik} which states that $\Ca$ must be equivalent to Rep $H$ for an algebraic supergroup $H$ whose representation category is semisimple. So we can assume $V$ to be a finite-dimensional (super) vector space. %Indeed, as $\Ca$ is generated by a single object, the closure in Zariski topology of the image of $H$ in $\End(V)$ generates the same tensor category as $H$.

If $H$ is a group, it follows from e.g. \cite{KLV}, Theorem 1.2 that $H$ must be the algebraic group $G$ of type $G_2$. Observe that a tensor category of type $G_2$  does not contain any finite subcategory  $\mathcal{F}$ except for $Vec$. Indeed, 
if  $\mathcal{F}$ contained an object $X_\la$ labeled by $\la\neq 0$, it would also contain simple objects $X_{n\la}$ for all $n\in{\mathbb N}$. Hence it only remains to consider the case when $H$ is a connected supergroup which is not a group. As Rep $H$ is semisimple, it is well-known that it has to be isomorphic to $SOSp(1|2n)$ for some $n$ (see e.g. \cite{CH}, Section 7 for a definition). It will be slightly easier to deal with
$OSp(1|2n)\cong SOSp(1|2n)\times \Z/2$. In \cite{CH} Theorem 7.19 it was shown that Rep $OSp(1|2n)$
is equivalent to the quotient of the Deligne category $O_{1-2n}$ modulo negligible morphisms.
By \cite{TW2}, Proposition 9.2, this quotient category has a fusion ring isomorphic to the one of Rep $O(2n+1)$. 
Hence Rep $H$ has a fusion ring isomorphic to the one of Rep $SO(2n+1)$. Again, it can not be isomorphic to Rep $G$, by \cite{KLV}, Theorem 1.2.
\subsection{Main result}\label{main:sec}  We now check that the assignment $F$ defined in Section \ref{functordef} is indeed a ribbon tensor functor between $\Ca$ and $\U_q$ resp $\bar{\U}_q$, see Section \ref{fusion:sec} for precise definitions.

\begin{theorem}\label{classification-theorem} 
Let $\Ca$ be a ribbon tensor category of type $G_2$ or $G_{2,k}$, $k\geq 1$, with $V$ the object corresponding to the 7-dimensional simple representation of the Lie algebra of type $G_2$. If $k=3$ or $6$, the classification is given by
Lemma \ref{lowtensor}. In the remaining cases, we denote the eigenvalue of $c_{V,V}$ via which it acts on the summand $V_{[2,0]}$ of $V^{\otimes 2}$ by $q^2$. Then we have

(a) If $\Ca$ is of type $G_2$, it is equivalent to the category $\U_q=$ Rep $U_q\g(G_2)$ for $q^2$ not a root of unity or $q^2=1$.

(b) If $\Ca$ is of type $G_{2,k}$, $k\geq 1$, it is equivalent to the quotient category $\bar{\U}_q$ of  Rep $U_q\g(G_2)$ for $q^2$ a primitive $(k+12)$-th root of unity.
\end{theorem}

$Proof.$ We will give a proof here if $\Ca$ is of type $G_2$, which we can assume to be not symmetric by Proposition \ref{sym:prop}. We have already defined assignments $X\in Ob(\Ca)\to F(X)\in Ob(\U_q)$ and $f\in \Hom(X,Y)\to F(f)\in \Hom(F(X),F(Y))$.
By Lemma \ref{Fwell-defined}, $F$ is a functor.
We now check that it preserves tensor products.

If $f:X_1\to Y_1$ and $g: X_2\to Y_2$ are morphisms, we observe that $i_{X_1}\otimes i_{X_2}$
defines an embedding of $X_1\otimes X_2$ into $V^{\otimes n_1}\otimes V^{\otimes n_2}$. Hence
$$F(f)\otimes F(g)=\dt_{Y_1}\Phi(i_{Y_1}fd_{X_1})\it_{X_1}\otimes \dt_{Y_2}\Phi(i_{Y_2}gd_{X_2})\it_{X_2}=$$ 
$$
(\dt_{Y_1}\otimes \dt_{Y_2})\Phi((i_{Y_1}\otimes i_{Y_2})(f\otimes g)(d_{X_1}\otimes d_{X_2}))\it_{X_1}\otimes \it_{X_2}=F(f\otimes g).$$
Our functor also maps our braiding morphism $c_{V,V}$ to the corresponding braiding morphism in $\U_q$ with the same eigenvalues. As the twisting numbers are already determined by the eigenvalues of $c_{V,V}$, see Proposition \ref{Casimirprop}, and the braiding structure is determined by the one of a generating object, our functor is a functor of ribbon categories.

The same proof works in type $G_{2,k}$, assuming uniqueness of the braid representation as it was proved in Theorem \ref{unique:thm} for type $G_2$. As discussed in Remark \ref{general case:rem}, this case will be proved in Section \ref{altunique} by a somewhat different method.

\subsection{Equivalences} We will assume Theorem \ref{classification-theorem} to be true also for type $G_{2,k}$ for which a complete proof will be given in Section \ref{altunique}. In order to get a complete classification of tensor categories of type $G_2$ or $G_{2,k}$, it remains to determine when the categories in Theorem \ref{classification-theorem} are equivalent. For simplicity, we refer to them just by $\Ca(q^2)$, with $q^2$ having one of the values allowed for the given fusion ring.
The following simple lemma uses an argument which already appeared before e.g. in \cite{KW}.

\begin{lemma}\label{noaut} Assume $\Ca(q^2)$ and $\Ca(\tilde q^2)$ are monoidally equivalent categories of type $G_2$ or $G_{2,k}$ with the same fusion ring, which does not allow any nontrivial automorphism. Then $\tilde q^2=q^{\pm 2}$.
\end{lemma}

$Proof.$ Let $p$ and $\tilde p$ be the projections onto
the summands labeled by $\La_2$ in $V^{\otimes2}$
in the categories $\Ca(q^2)$ and $\Ca(\tilde q^2)$.
A monoidal equivalence preserves this summand,
since the fusion ring has no nontrivial automorphism. Let $p_1=p\otimes 1$ and $p_2=1\otimes p$.
It follows from Proposition \ref{TLnew} that $p_1p_2p_1=\frac{1}{[2]^2}p_1$ in $\Hom(V_{[2,1]},V^{\otimes 3})$. Then an equivalence between $\Ca(q^2)$ and $\Ca(\tilde q^2)$ would map $p_i$ to $\tilde p_i$, for $i=1,2$. This implies $1/(q+q^{-1})^2=1/(\tilde q+\tilde q^{-1})^2$. It is straightforward to check that this is possible only if $\tilde q^2=q^{\pm 2}$. As $\Ca(q^{-2})$ is just $\Ca(q^2)$ with the opposite braiding, these two categories are indeed equivalent as monoidal categories.

\begin{lemma}\label{G29} Let $\Ca$ be a tensor category of type $G_{2,9}$, and let $q^2$ be a primitive 21st root of unity. Then $\Ca(q^{\pm 2})\sim \Ca(q^{\pm 8})\sim \Ca(q^{\pm 10})$.
In particular, there are only two non-equivalent monoidal tensor categories of type $G_{2,9}$.
\end{lemma}

$Proof.$ It is known, see \cite{Ga} that the fusion ring $G_{2,9}$ has a nontrivial automorphism $\alpha$ of order three given by
$$[1,0]\ \mapsto \ [1,1]\ \mapsto\ [3,0]\  \mapsto\ [1,0],$$
where the diagrams [0,0], [2,0] and [2,1] are fixed. We claim that this automorphism induces a functor from $\Ca(q^2)$ to $\Ca(q^{-10})$. 
It follows from our fusion ring automorphism for the object $W=V_{[1,1]}$ that
$W\otimes W\ \cong \1\oplus W\oplus V_{[3,0]}\oplus V_{[2,0]}.$
Observe that the eigenvalue of $c_{W,W}^2$ on $V_{[2,0]}$ is given by $q^{c_{[2,0]}-2c_W}=q^{28-48}=q^{-20}$, see Eq. \ref{ribbon} and Proposition \ref{prop:eigenvalue}. Hence $c_{W,W}$ acts on $V_{[2,0]}$ via $q^{-10}$; it cannot act by $-q^{-10}$, as $W$ generating a category of type $G_{2,9}$ requires the eigenvalue for $V_{[2,0]}$ to be a primitive 21st root of unity. This shows that $\Ca(q^{-10})=\Ca(q^{-5\cdot 2})$ is equivalent to $\Ca(q^2)$. 
We show that $\Ca(q^{-5(-10)})=\Ca(q^8)$ is equivalent to $\Ca(q^{-10})$ by the same argument. By \cite{Ga}, Theorem 3.G2, the only fusion ring
automorphisms are $1,\alpha,\alpha^2$.
The equivalences above already give these changes of labels. Given any other equivalence, we can compose it with the inverse of one of these so that it fixes every simple label. The calculation in the proof of Lemma \ref{noaut} then shows that the parameters must be equal or inverse. Thus the displayed equivalences give all possible
identifications. The twelve primitive 21st roots
form two orbits of six, giving exactly two monoidal equivalence classes.

\begin{theorem}\label{explicit classification} Let $\Ca(q^2)$ denote the (quotient) categories of $\U_q$ in Theorem \ref{classification-theorem}. Then $\Ca(q^2)$ is equivalent to $\Ca(\tilde q^2)$ as a monoidal tensor category if and only if $\tilde q^2=q^{\pm 2}$, except for the case stated in Lemma \ref{G29}.
\end{theorem}

$Proof$. By \cite{Ga}, the fusion ring $G_{2,k}$ has nontrivial automorphisms only for $k=9$ and $k=12$. We have already dealt with $k=9$ in Lemma \ref{G29}. It was shown in \cite{EM}, Lemma 6.4 that the automorphism for $k=12$ induces a braided autoequivalence of $\Ca(q^2)$. Hence we get no new identifications.

It is known that the fusion ring for $G_2$ has no nontrivial automorphisms. This can e.g. be seen as follows. It is well-known that for any object $W$ in $\Ca$ the tensor product $V_\la\otimes W$ decomposes into the direct sum of $\dim W$ simple representations if the dominant integral weight $\la$ is far enough from the walls of the dominant Weyl chamber; here $\dim W$ is the classical dimension for the Lie algebra $G_2$. Hence a fusion ring automorphism must map a simple object to an object with the same classical dimension. It follows from the dimension formulas \ref{dimensionformula} that the only objects with dimensions 7 and 14 are the two fundamental representations. The claim hence follows from this and Lemma \ref{noaut}.
\subsection{Trivalent categories} A graphical description of the tensor category $\U_q$ of representations of the Drinfeld-Jimbo quantum group of type $G_2$ was given by Kuperberg \cite{Ku2}, and further studied in \cite{MPS} and \cite{MST}. It exploits the existence of a nontrivial morphism $\tau: V\to V^{\otimes 2}$ in such categories. We give a brief description of it here, based on the presentation in \cite{MST}, Definition 2.15. Moreover, we make precise how one can get the description in terms of trivalent graphs from our braid representations. We use this to give an alternative and complete proof of Theorem \ref{classification-theorem}.

We start out with a preliminary category $\tilde\G(q)$ whose objects are the non-negative integers.
The morphisms are given by graphs generated by the usual braiding graph $\st$ (for a chosen over- or undercrossing) and duality graphs $\cup$ and $\cap$, as well as trivalent graphs $\tt$ and $\tt^t$, see below. Then
$\Hom(m,n)$ is given by all linear combinations of such graphs which connect $m$ lower vertices with $n$ upper vertices, subject to a number of relations, see 
\cite{MST}, Definition 2.15. Several of these relations will also appear below.
As usual, the concatenation of morphisms $\ft\circ \gt$ is given by stacking the graph $\ft$ on top of the graph $\gt$.  The category $\G(q)$ is  defined as the Karoubian envelope, whose objects are idempotents in $\End(n)$ for all $n\in \N$. We then have

\begin{theorem}\label{Kuperbergtheorem} \cite{Ku2} Let $V$ be the object corresponding to the simple 7-dimensional object in $\U_q$ for $q$ not a root of unity except $\pm 1$. Then we have linear isomorphisms between $\Hom(m,n)$ and $\Hom(V^{\otimes m},V^{\otimes n})$
which define an equivalence between $\G(q)$ and $\U_q$ for $q^2$ not a root of unity, except for 1.
\end{theorem}
The linear maps in this theorem are determined by mapping
the braiding and duality graphs $\cup$ and $\cap$ in $\G(q)$ to the braiding morphism $c_{V,V}$ and the duality morphisms $\1\to V^{\otimes 2}$ and $V^{\otimes 2}\to \1$ in $\U_q$ for the self-dual object $V$, with the normalization 
$$\cap\circ\cup=d=q^{10}+q^8+q^2+1+q^{-2}+q^{-8}+q^{-10}.$$
The assignments for the trivalent graphs are listed below, where $\tau$ and $\tau^t$ are nonzero morphisms normalized such that  $\tau^t\circ\tau=b$, with $b=-[3](q^4+q^{-4})$   as in \cite{MST}.
\[ \tt^t\ =\ 
\begin{tikzpicture}[scale=0.6,baseline=-0.5ex,line cap=round,line join=round]
  % upside-down Y : V⊗V → V
  \draw[line width=0.8mm]
    (-1,-1) -- (0,0) -- (1,-1)
    (0,0) -- (0,1.5);
\end{tikzpicture}
\ \to \  \tau^t:\  V \otimes V\to V, \hskip 3em \tt\ =\ 
\begin{tikzpicture}[scale=0.6,baseline=-0.5ex,line cap=round,line join=round]
  % upright Y : V → V⊗V
  \draw[line width=0.8mm]
    (0,-1.5) -- (0,0) -- (-1,1)
    (0,0) -- (1,1);
\end{tikzpicture}
\ \to \  \tau:\  V \to V\otimes V.
\]
We also have the following morphism in the description of the relations for $\G(q)$.
\[
\at \;=\;
\begin{tikzpicture}[scale=0.5,baseline=-0.5ex,line cap=round,line join=round]
  \begin{scope}[rotate=90]
    \draw[line width=0.9mm]
      (-1,1) -- (0,0) -- (1,1)
      (0,0) -- (0,-1.2) -- (-1,-2.2)
      (0,-1.2) -- (1,-2.2);
  \end{scope}
\end{tikzpicture}
\; = \;\begin{tikzpicture}[scale=0.5,baseline=-0.5ex,
    line cap=round,line join=round]

  % =====================
  % Top layer 
  % baseline y = 1
  \draw[line width=0.9mm] (-3,1) -- (-3,1.8);
  \draw[line width=0.9mm] (-1,1) -- (-1,1.8);
  \draw[line width=0.9mm] (1,1) arc[start angle=180,end angle=0,radius=1];

  % =====================
  % Middle layer (b with vertical sticks)
  \draw[line width=0.9mm] (-3,-2.2) -- (-3,1);  % left stick continuous
  \draw[line width=0.9mm] ( 3,-2.2) -- ( 3,1);  % right stick continuous
  % the b
  \draw[line width=0.9mm]
    (-1,1) -- (0,0) -- (1,1)
    (0,0) -- (0,-1.2) -- (-1,-2.2)
    (0,-1.2) -- (1,-2.2);

  % =====================
  % Bottom layer, baseline y = -2.2
  \def\Y{-2.2} \def\R{1.0}
  \draw[line width=0.9mm] (-3,\Y) arc[start angle=180,end angle=360,radius=\R];
  \draw[line width=0.9mm] ( 1,\Y-\R) -- ( 1,\Y);
  \draw[line width=0.9mm] ( 3,\Y-\R) -- ( 3,\Y);

\end{tikzpicture}
\]
We can express $\tt\circ\tt^t$ and $\cup\circ\cap$ as multiples of eigenprojections of the braiding graph $\st$. We can also express the graph $\at$ just using our braid representations as
\begin{equation}\label{aexpression}
\a\otimes (\cup\circ \cap)\ =\ (1_2\otimes (\cup\circ \cap))\ \circ\ (1\otimes \ (\tt\circ\tt^t) \otimes 1)\ \circ\ ((\cup\circ \cap)\otimes 1_2)\ \circ\ \st_2\st_1\st_3\st_2,
\end{equation}
where $\st_i$ is the graph for the $i$-th braid generator $\sigma_i$.

\def\Cat{\tilde{\mathcal C}}
\subsection{Trivalent relations from braids} We now relate categories of type $G_2$ or $G_{2,k}$ to the categories $\G(q)$, using Lemma \ref{noah}.

\begin{lemma}\label{trivalent functor} Let $\Ca$ be a tensor category of type $G_2$ or $G_{2,k}$, with $q^2\neq\pm1$ the eigenvalue of $c_{V,V}$ corresponding to the object labeled by $2\La_1=[2,0]$. Then there is a functor $G:\G(q)\to \Ca$,
which is already uniquely determined by the representations $\C B_j\to \End(V^{\otimes j})$, $j\leq 4$.
\end{lemma}

$Proof.$ The idea is to express the generating graph morphisms in the framework of our braid representations and to check that these satisfy the defining relations. This can be done within $B_4$. By our assumptions on $\Ca$, there exist morphisms $\tau: V\to V^{\otimes 2}$ and $\tau^t: V^{\otimes 2}\to V$, which we normalize such that
$\tau^t\circ \tau=b$. We define the functor $G$ by the same assignments as in the previous subsection, now for $\Ca$ instead of $\U_q$.
It remains to check that these assignments preserve the relations for $\G(q)$.  We are going to show this for the two most complicated relations. One of them is given by the picture below, where $a=q^4+1+q^{-4}$.

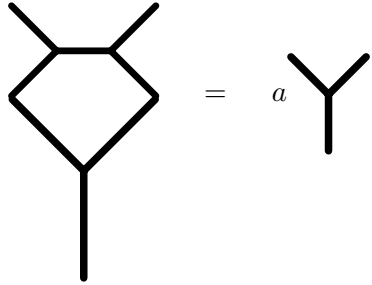
\[
\begin{tikzpicture}[baseline=-4ex,line cap=round,line join=round]

  % --- TOP: ---
  \begin{scope}[scale=0.6]
    \begin{scope}[rotate=90]
      \draw[line width=0.9mm]
        (-1,1) -- (0,0) -- (1,1)
        (0,0) -- (0,-1.2) -- (-1,-2.2)
        (0,-1.2) -- (1,-2.2);
    \end{scope}
  \end{scope}

  % --- BOTTOM: Y ---
  \begin{scope}[scale=0.95,yshift=-1.68cm, xshift=0.38cm]
    \draw[line width=1mm]
      (0,-1.5) -- (0,0) -- (-1,1)
      (0,0) -- (1,1);
  \end{scope}
\end{tikzpicture}
\hskip 1em \;=\;\hskip 1em 
a \begin{tikzpicture}[scale=0.5,baseline=-0.5ex,line cap=round,line join=round]
  % Y
  \draw[line width=1mm]
    (0,-1.5) -- (0,0) -- (-1,1)
    (0,0) -- (1,1);
\end{tikzpicture}
\]
Let $G(\at)=\al\in \End(V^{\otimes 2})$. Then we have to check $\alpha\circ\tau= a \tau$.
This just means that the eigenvalue of $\alpha$ corresponding to the summand $V\subset V^{\otimes 2}$ is equal to $a$. This is true for $\Ca=\U_q$, by Theorem \ref{Kuperbergtheorem} for $q$ not a root of unity. It extends to roots of unity in an obvious way by continuity. But as $G$ was defined in terms of braid representations for $\Ca$  which are isomorphic to the ones for $\U_q$ by our assumptions, it is also true for $\Ca$. 

To formulate the second relation, let  $M_2=\{\at, \tt\circ\tt^t, \cup\circ\cap, 1_2\}$. 
Let $\sqt$ in $\End(2)$ be $\frac{1}{d}$ times the graph given by a square with trivalent edges at the corners (see \cite{MST}, Definition 2.15 or the formula below). Then the relation states that
\begin{equation}\label{sqrel}
\sqt=\sum_{\bt\in M_2} a_\bt \bt
\end{equation}
for certain scalars $a_\bt$, i.e. $\sqt$ is a linear combination of the elements of $M_2$.
As $\pb_o=\frac{1}{d}\cup\circ\cap$ is the projection onto the summand in the object $2\in Ob \G(q)$ which is isomorphic to the trivial object, we can embed $\End(2)$ into $\End(4)$ via
the map
$$\bt\in \End(2)\ \mapsto\ \hat\bt\ = \ \frac{1}{d}(1\otimes \cup\ \otimes 1) \ \circ\ \bt\  \circ\ (1\otimes \cap\ \otimes 1)\ =\ \st_2\st_3\ \circ\ (\bt\otimes \pb_o)\ \circ\  \st_3^{-1}\st_2^{-1};$$
informally,  $\hat\bt$ is $\bt$ acting on the first and last vertex in $\End(4)$, and $\pb_o$ acts on the two middle ones.
One similarly expresses (see the defining picture of $\sqt$ in \cite{MST}, Definition 2.15)
$$\hat\sqt\ =\ \ (1\otimes \pb_o\otimes 1)\ \circ\ (\tt\circ\tt^t\otimes \tt\circ\tt^t)\ \circ\  (1\otimes \pb_o\otimes 1).$$
Using our embedding $\bt\in \End(2)\to \hat\bt\in \End(4)$, Relation \ref{sqrel} becomes $\hat\sqt=\sum_{\bt\in M_2} a_\bt \hat\bt$. It holds in $\U_q$. But as the images of $\hat\sqt$ and $\hat\bt$ in $\Ca$ only depend on the representation of $B_4$ in $\End(V^{\otimes 4})$, it also must hold in $\Ca$.
The other relations for $\G(q)$ are much easier. The reader should have no problem checking them for themselves.

\subsection{Complete proof of Theorem \ref{classification-theorem}}\label{altunique} We can now give a proof of Theorem \ref{classification-theorem} which also covers categories of type $G_{2,k}$. By Lemma \ref{noah}, we have uniqueness of braid representations into $\End(V^{\otimes j})$ for $j\leq 4$ for any tensor category of type $G_2$ or $G_{2,k}$ provided $\la_1=q^2\neq\pm1$ is fixed. But then Lemma \ref{trivalent functor} shows that we obtain a functor of $\G(q)$ into $\Ca$. As $\G(q)$ contains the braiding morphisms, the morphisms from $\G(q)$ already generate $\End(V^{\otimes n})$ for all $n\in \N$, by Theorem \ref{surjective:thm}. Moreover, as $\Ca$ does not contain any negligible morphisms, it follows that $\G(q)/\Neg$ maps onto $\Ca$; here $\Neg$ is the tensor ideal of negligible morphisms in $\G(q)$. As $\Neg$ is a maximal tensor ideal (see \cite{AK}, or \cite{MPS}, Proposition 2.3), it follows that $\Ca$ is equivalent to  $ \G(q)/\Neg$. Hence $q^2$ determines the tensor category of type $G_2$ or $G_{2,k}$ up to equivalence, which proves Theorem 
\ref{classification-theorem}. 

\medskip

As an immediate consequence of our proof, we can rephrase our main result in terms of the category $\G(q)$. In particular, this answers one of the questions discussed in \cite{MPS}, Remark 5.2.6. for the category $\G(q)$ for $q$ a root of unity.

\begin{theorem}\label{quotient:corollary} Let $\Ca$ be a tensor category of type $G_2$ or $G_{2,k}$ with $k\geq 1$ and $\la_1=q^2\neq \pm 1$.
Then $\Ca$ is equivalent to $\G(q)/\Neg$, with $q$ as determined in Theorem \ref{dimfromeigen}.
\end{theorem}
\subsection{Discussion of results and future directions}\label{discussion} 1. Our classification of tensor categories was based on certain rigidity properties of braid representations. We first use the classification of simple representations of $B_4$ for which the generators have exactly three eigenvalues (see \cite{MW3}) to explicitly determine the eigenvalues for the braiding morphism $c_{V,V}$, see Proposition \ref{eigenvaluesribbon} and Theorem \ref{dimfromeigen}.
Using path representations and a generalization of the quantum Jucys-Murphy approach,
see Theorem \ref{AB2theorem}, we can then show that the image of the braid group $B_n$ generates $\End(V^{\otimes n})$ for all $n\in \N$ as long as the category $\Ca$ of type $G_2$ or $G_{2,k}$ is not symmetric. One then derives our main result Theorem \ref{classification-theorem} either by proving that these braid representations are unique, see Theorem \ref{unique:thm} and Remark \ref{general case:rem}, or by relating them to Kuperberg's spiders, see Lemma \ref{noah} and Section \ref{altunique}. If $\Ca$ is a symmetric tensor category of type $G_2$, it has been known to experts that it is unique up to equivalence, see Section \ref{symmetric:sec}.

2.  The possibility of classifying tensor categories of type $G_2$ was already considered by Morrison, Snyder and Peters, see the remark before Corollary 8.9 in \cite{MPS}. Their approach, which goes back to Kuperberg's spiders in \cite{Ku2} is quite different from ours. Our approach is primarily based on analyzing braid representations, while theirs depends on the existence of a nontrivial morphism $\tau: V\to V^{\otimes 2}$. It is possible to prove our main result in many cases without using Kuperberg's formalism, see Remark \ref{general case:rem}. We did use it to give a complete proof of Theorem \ref{classification-theorem} without excessive additional calculations. But, at this point, we do not see a short way to prove Theorem \ref{surjective:thm} at roots of unity with only using the spider approach.  It implies, in particular, that  categories of type $G_2$ or $G_{2,k}$ with $k\geq 1$ are equivalent to quotients of a trivalent category $\G(q)$ modulo its maximal ideal $\Neg$ of negligible morphisms, see Theorem \ref{quotient:corollary}. 
There are also interesting quotient categories of $\G(q)$ for $q^2$ a primitive fourth, fifth, and ninth root of unity. We will study them in \cite{RW}.

3. We would like to thank Andrew Schopieray for pointing out that our classification of tensor categories of type $G_{2,9}$  solves one of the open problems left in the classification of ribbon fusion tensor categories of rank 6 in \cite{DS}.

4. Tensor categories of the classical Lie types $A_n$, $C_n$, $B_n$ and $D_n$ have been classified in \cite{FK}, \cite{KW}, \cite{TW} and \cite{Cp}, where a ribbon structure is assumed except for type $A$.
This does not include the spinor representations for types $B$ and $D$ so far. We expect that the results in this paper and in \cite{MW3} should also be useful for studying other exceptional Lie types. E.g. if $V$ is the adjoint representation, or for $E_6$ and $E_7$ the smallest nontrivial representation,
the images of the generator $\sigma_i$ in the braid representations into $\End(V^{\otimes n}_{new})$ only have three eigenvalues. In the latter two cases, these representations were already studied in \cite{Wexc}. Although other exceptional types will be more complicated (e.g., the braid groups do not even generate $\End(V^{\otimes n})_{new}$ in general), a combination of results on braid representations in this paper, \cite{MW3} and \cite{MW2} with recent results in \cite{MST} should help to get a better understanding of such tensor categories. In the latter paper, the authors Morrison, Snyder and Thurston consider a possible quantum version of Deligne's conjectured exceptional series, see \cite{D2}. 
Regardless whether this series exists or not, their approach is useful for giving a graphical description of tensor categories of exceptional types,  similarly as the special case of Kuperberg's spiders is useful for $G_2$.

\end{document}